\input amstex
\documentstyle{amsppt}
\magnification =\magstep1
\baselineskip=18pt
\vcorrection{-.33truein}
\pageheight{9.0truein}
\input EPSfig.macros
\nologo
\overfullrule=0pt

\newif\ifproofmode                      
\proofmodefalse				

\newif\ifforwardreference		
\forwardreferencefalse			

\newif\ifchapternumbers			
\chapternumbersfalse			

\newif\ifcontinuousnumbering		
\continuousnumberingfalse		

\newif\iffigurechapternumbers		
\figurechapternumbersfalse		

\newif\ifcontinuousfigurenumbering	
\continuousfigurenumberingfalse		

\font\eqsixrm=cmr6			
\def\marginstyle{\eqsixrm}		

\newtoks\chapletter			
\newcount\chapno			
\newcount\eqlabelno			
\newcount\figureno			

\chapno=0
\eqlabelno=0
\figureno=0


\def\chapfolio{\ifnum \chapno>0 \the\chapno \else \the\chapletter \fi}


\def\bumpchapno{\ifnum \chapno>-1 \global \advance \chapno by 1
	\else \global \advance \chapno by -1 \setletter\chapno \fi
	\ifcontinuousnumbering \else \global\eqlabelno=0 \fi
	\ifcontinuousfigurenumbering \else \global\figureno=0 \fi}

%


%

\def\tempsetletter#1{\ifcase-#1 {}\or{} \or\chapletter={A}\or\chapletter={B}
  \or\chapletter={C} \or\chapletter={D} \or\chapletter={E}
  \or\chapletter={F} \or\chapletter={G} \or\chapletter={H}
  \or\chapletter={I} \or\chapletter={J} \or\chapletter={K}
  \or\chapletter={L} \or\chapletter={M} \or\chapletter={N}
  \or\chapletter={O} \or\chapletter={P} \or\chapletter={Q}
  \or\chapletter={R} \or\chapletter={S} \or\chapletter={T}
  \or\chapletter={U} \or\chapletter={V} \or\chapletter={W}
  \or\chapletter={X} \or\chapletter={Y} \or\chapletter={Z}\fi}

%

\def\chapshow#1{\ifnum #1>0 \relax #1%
   \else {\tempsetletter{\number#1}\chapno=#1 \chapfolio} \fi}

%
\def\today{\ifcase\month\or
January\or February\or March\or April\or May\or June\or
July\or August\or September\or October\or November\or December\fi
\space\number\day, \number\year}

%

\def\initialeqmacro{\ifproofmode
 \headline{\tenrm \today\hfill \jobname\ --- draft\hfill\folio}
     \hoffset=-1cm \immediate\openout2=allcrossreferfile \fi
 \ifforwardreference \input labelfile
     \ifproofmode \immediate\openout1=labelfile \fi \fi}


%

\def\chaplabel#1{\bumpchapno \ifproofmode \ifforwardreference
   \immediate\write1{\noexpand\expandafter\noexpand\def
   \noexpand\csname CHAPLABEL#1\endcsname{\the\chapno}}\fi\fi
   \global\expandafter\edef\csname CHAPLABEL#1\endcsname
   {\the\chapno}\ifproofmode\llap{\hbox{\marginstyle #1\ }}\fi\chapfolio}

%
\def\eqnum{\global\advance\eqlabelno by 1
   \eqno(\ifchapternumbers\chapfolio.\fi\the\eqlabelno)}

\def\eqlabel#1{\global\advance\eqlabelno by 1 \ifproofmode\ifforwardreference
 \immediate\write1{\noexpand\expandafter\noexpand\def
 \noexpand\csname EQLABEL#1\endcsname{\the\chapno.\the\eqlabelno?}}\fi\fi
 \global\expandafter\edef\csname EQLABEL#1\endcsname
 {\the\chapno.\the\eqlabelno?} \eqno(\ifchapternumbers\chapfolio.\fi
 \the\eqlabelno)\ifproofmode\rlap{\hbox{\marginstyle #1}}\fi}

\def\leqlabel#1{\global\advance\eqlabelno by 1 \ifproofmode\ifforwardreference
 \immediate\write1{\noexpand\expandafter\noexpand\def
 \noexpand\csname EQLABEL#1\endcsname{\the\chapno.\the\eqlabelno?}}\fi\fi
 \global\expandafter\edef\csname EQLABEL#1\endcsname
 {\the\chapno.\the\eqlabelno?} \leqno(\ifchapternumbers\chapfolio.\fi
 \the\eqlabelno)\ifproofmode\rlap{\hbox{\marginstyle #1}}\fi}

\def\eqalignnum{\global\advance\eqlabelno by 1
   &(\ifchapternumbers\chapfolio.\fi\the\eqlabelno)}

\def\eqalignlabel#1{\global\advance\eqlabelno by1 \ifproofmode
 \ifforwardreference\immediate\write1{\noexpand\expandafter\noexpand\def
 \noexpand\csname EQLABEL#1\endcsname
     {\the\chapno.\the\eqlabelno?}}\fi\fi
 \global\expandafter\edef\csname EQLABEL#1\endcsname
 {\the\chapno.\the\eqlabelno?}\ifchapternumbers\chapfolio.\fi
 \the\eqlabelno\ifproofmode\rlap{\hbox{\marginstyle
 #1}}\fi}

\def\eqref#1{(\ifundefined{EQLABEL#1}***\ifproofmode\ifforwardreference)%
   \else \write16{ ***Undefined Equation Reference #1*** }\fi
   \else \write16{ ***Undefined Equation Reference #1*** }\fi
   \else \edef\LABxx{\getlabel{EQLABEL#1}}%
   \def\LAByy{\expandafter\stripchap\LABxx}\ifchapternumbers
   \chapshow{\LAByy}.\expandafter\stripeq\LABxx
   \else\ifnum \number\LAByy=\chapno \relax\expandafter\stripeq\LABxx
   \else\chapshow{\LAByy}.\expandafter\stripeq\LABxx\fi\fi)\fi
   \ifproofmode\write2{Equation #1}\fi}

%

\def\fignum{\global\advance\figureno by 1 \relax
   \iffigurechapternumbers\chapfolio.\fi\the\figureno}\

\def\figlabel#1{\global\advance\figureno by 1\relax
 \ifproofmode\ifforwardreference
 \immediate\write1{\noexpand\expandafter\noexpand\def
 \noexpand\csname FIGLABEL#1\endcsname{\the\chapno.\the\figureno?}}\fi\fi
 \global\expandafter\edef\csname FIGLABEL#1\endcsname
 {\the\chapno.\the\figureno?}\iffigurechapternumbers\chapfolio.\fi
 \ifproofmode$^{\hbox{\marginstyle #1}}$\relax\fi\the\figureno}

\def\figref#1{\ifundefined{FIGLABEL#1}!!!!\ifproofmode\ifforwardreference)%
   \else \write16{ ***Undefined Equation Reference #1*** }\fi
   \else \write16{ ***Undefined Equation Reference #1*** }\fi
   \else \edef\LABxx{\getlabel{FIGLABEL#1}}%
   \def\LAByy{\expandafter\stripchap\LABxx}%
   \iffigurechapternumbers\chapshow{\LAByy}.\expandafter\stripeq\LABxx
   \else\ifnum\number\LAByy=\chapno \relax\expandafter\stripeq\LABxx
   \else\chapshow{\LAByy}.\expandafter\stripeq\LABxx\fi\fi
   \ifproofmode\write2{Figure #1}\fi\fi}

%

%

\def\getlabel#1{\csname#1\endcsname}
\def\ifundefined#1{\expandafter\ifx\csname#1\endcsname\relax}
\def\stripchap#1.#2?{#1}
\def\stripeq#1.#2?{#2}

\figurechapternumberstrue  

\chapternumberstrue        

\def\thmlbl#1{\figlabel{#1}}
\def\thmref{\figref}
\def\eqnlbl#1{\leqlabel{#1}}

\def\eqnref#1{\eqref{#1}}
\def\sectionnumber{\chapno}
\def\theoremnumber{\figureno}
\def\equationnumber{\eqlabelno}

\input amssym.def
\def\notexist{\raise1pt\hbox{$\not$}\exists}

\def\BbbR{\Bbb{R}}



\def\sgn{\text{\rm \ sgn \ }}
\def\Tr{\text{\rm \ Tr \ }}

\def\myqed{\vrule height3pt depth2pt width3pt \bigskip}
\def\newsection{\centerline}
\def\const.{\roman{const.}}

\def\R{\roman{Re\ }}



\def\barU{{\bar{u}}}

\def\Zbar{{\bar{z}}}

\def\T3{{\widetilde{\xi}}}

\def\Tsigma{{\widetilde{\sigma}}}

\def\tildeX{{\widetilde{x}}}

\def\tildeY{\widetilde{y}}

\pagewidth{6truein}
\pageheight{9truein}
\topmatter
\title
A stability index for
detonation waves in 
Majda's model for reacting flow
\endtitle
\leftheadtext{A stability index for detonation waves}
\rightheadtext{Gregory Lyng and Kevin Zumbrun}
\thanks
Research of both authors was supported
in part by the National Science Foundation under Grants No. DMS-9107990
and DMS-0070765.
\endthanks

\abstract
Using Evans function techniques,
we develop a stability index for 
weak and strong detonation waves
analogous to that developed for shock waves in [GZ,BSZ],
yielding useful necessary conditions for stability.
Here, we carry out the analysis in the context of the Majda model,
a simplified model for reacting flow;
the method is extended to the full Navier--Stokes equations of
reacting flow in [Ly,LyZ].
The resulting stability condition is satisfied for all nondegenerate,
i.e., spatially exponentially decaying, weak and strong detonations 
of the Majda model  in agreement with numerical experiments of [CMR] and 
analytical results of [Sz,LY] for a related model of Majda and Rosales.
We discuss also the role in the ZND limit of degenerate, 
subalgebraically decaying weak detonation and (for a modified,
``bump-type'' ignition function) deflagration profiles,
as discussed in [GS.1--2] for the full equations. 
\endabstract
\author
{Gregory Lyng and Kevin Zumbrun}
\endauthor
\date{October 6, 2000; Revised: May 28, 2003}
\enddate
\address
Department of Mathematics
University of Michigan
Ann Arbor, MI 48109--1109
\endaddress
\email
glyng\@umich.edu
\endemail 
\address
Department of Mathematics,
Indiana University,
Bloomington, IN  47405-4301
\endaddress
\email
kzumbrun\@indiana.edu
\endemail 
\endtopmatter
\document
\newsection {\bf Section 1. Introduction}
\sectionnumber=1 
\theoremnumber=0
\equationnumber=0
\smallskip
\TagsOnLeft

\bigskip
In one-dimensional, Lagrangian coordinates, the 
Navier--Stokes equations of reacting flow for a
one-step reaction
may be written in the abstract form
$$
\cases
u_t + f(u)_x = (B(u)u_{x})_x +kq\varphi(u)z,\\
z_t = 
(D(u,z)z_x)_x
-k\varphi(u)z,
\endcases \eqnlbl{1abstract}
$$
where $u$, $f$, $q\in \BbbR^n$, $B\in \BbbR^{n\times n}$,
$z$, $k$, $D$, $\varphi \in \BbbR^1$, and $k>0$
(model (8.81) of [Z.3] 
with particle velocity set to zero).
Here, vector $u$ comprises the gas-dynamical variables of
specific volume, particle velocity, and total energy, and $z$
measures mass fraction of unburned reactant: more generally,
``progress'' of a single reaction involving multiple reactants.
The first equation thus models
kinematic and the second equation reaction effects.  
The function $\varphi(u)$ is an ``ignition function'', monotone
increasing in temperature, and usually
assumed for fixed density to be zero below a certain ignition temperature 
and positive above.
The vector $q$ comprises quantities produced in reaction, 
in particular heat released. 
The coefficient $k$ corresponds to reaction rate, while
coefficients $B$ and $D$ model transport effects of, respectively,
viscosity and heat conduction, and species diffusion.
Multi-step reactions may be modeled by the same equations
with vectorial reaction variable $z\in \BbbR^m$, 
and coefficients $q$, $D$, $\varphi$, $k$
modified accordingly; see Section 9.
For further discussion, see, e.g.,  [CF,FD,GS.1--2,Z.3,Ly].

Under different conditions at $x=\pm
\infty$, there can result a variety of types of 
waves solving \eqnref{1abstract}:
nonreactive gas-dynamical shock and rarefaction solutions
($z\equiv 0$, or $z\equiv constant$ and $\varphi\equiv 0$), 
and traveling combustion waves consisting
of weak and strong detonations, weak and strong deflagrations, and
Chapman--Jouget detonation and deflagration waves, which are
limiting cases dividing weak and strong branches.  
Roughly speaking, detonations are compressive waves analogous to shock
waves in nonreactive gas dynamics, while deflagrations are expansive
solutions analogous to rarefactions;  for a fixed left-hand state,
there are weak and strong branches of right-hand sides corresponding
to waves of each type.  Chapman-Jouget waves occur at the special
parameters for which strong and weak branches coalesce.
We refer the reader to [CF,FD,G,M.4,GS.1--2]
for a detailed discussion of these solutions of the traveling wave
ODE for \eqnref{1abstract} and their roles in Riemann solutions/time-asymptotic
behavior for the initial value problem
under various assumptions on $B$, $D$.  

Similarly as in the case of ``real'', e.g., van der Waals, gas dynamics
[BE,MeP], the multitude of possible such elementary waves
leads to a multitude of possible time-asymptotic states, and these
must be classified according to stability.
However, the assessment of stability is a complicated undertaking.
Up to now, essentially all analyses have been carried out for one of
three simplified models: 
(i) the Zeldovich--von Neumann-Doering 
(ZND) model, for which $B$ and $D$ are set identically zero 
in \eqnref{1abstract} [Er.1--6,CF,LS], 
(ii) the still further simplified Chapman--Jouget (CJ)
or ``square-wave'' model, for which 
$B$ and $D$ are set identically zero and $k$ is taken to be infinite, 
so that combustion waves become surfaces of discontinuity 
[Er.7,F,etc.],
or (iii) The Majda model, for which $u$ is taken to be a scalar, and
$B$ and $D$ are set to $1$ and $0$ (or sometimes $1$), respectively
[M.4,LLT,LYi,RV],
or the Majda--Rosales model [CMR,Sz,LY,Li.1--7],
a closely-related cousin
in which $-z_x$ is substituted for $z_t$ everywhere in \eqnref{1abstract}.
%
At one time, it seems to have been believed,
based on analysis of the ZND case
$B\equiv D\equiv 0$ (see, e.g., [CF])
that weak detonations and strong deflagrations were unstable, the other
types stable at least in moderate parameter ranges.  
This conjecture on weak detonation is now widely
agreed to be false when viscosity and other effects
are taken into account, see [CMR,Sz,LY], 
or more general discussion in [FD].  
However, rigorous analysis of stability for the full model \eqnref{1abstract},
or comparison with stability for the ZND or CJ approximations, remain
important open problems.

The purpose of the present paper is to 
initiate a larger-scale study of these problems by the introduction,
in the simple setting of the Majda model, of 
new Evans function techniques 
developed recently
in the study 
of stability of viscous shock profiles
(see, e.g., [GZ,ZH,BSZ,ZS,Z.3]).
In contrast to the methods of past analyses of the Majda model,
these techniques were designed for the study of {\it systems},
$u\in \BbbR^n$, $n>1$, 
so may be applied also in the case of the
full, reacting Navier--Stokes equations.
As a roadmap for discussions to follow,
we point out that the ZND limit 
$B=\varepsilon B_0$, $D=\varepsilon D_0$,  $\varepsilon\to 0$ 
is equivalent by the spatial rescaling $x\to x/\varepsilon$ 
to the small-$k$ limit $k\to 0$ with $B$, $D$ held fixed.
The CJ, or square-wave limit in the viscous setting 
($B$, $D$ fixed) is ambiguous,
corresponding to intermediate values of $k$, but large spatial scale;
see discussion of [Z.3], Appendix A.3.
The large $k$ limit $k\to \infty$ (with, necessarily, $u_i$, $u^i\to u_-$
to allow a connection) corresponds roughly to
the high activation-energy limit described in, e.g., [B], pp. 22--25, 
which in the ZND setting leads to a square-wave 
approximation theory, and is generally associated with instability, 
and other anomalous behavior; see, e.g., [Er.7,AT,BL,BN,LS].

A tool that has proved useful for the study of stability
in the related cases of van der
Waals gas dynamics and multiphase flow [GZ,Z.6]
is a one-dimensional {\it stability index} originally
introduced by J. Evans [E.1--4] in the context of nerve axon equations, 
and generalized in various directions in, e.g., [J.1,AGJ,PW,GZ,BSZ].
A topological index relating evolutionary (PDE) dynamics of a 
traveling wave to dynamics of the associated traveling wave ODE,
the stability index is based on the {\it Evans function} 
[E.1--4,AGJ,PW,GZ,BSZ,etc.] $ D(\lambda)$, 
an analytic function playing the role of a characteristic function
for the linearized operator $L$ about the wave.
Defined as a Wronskian of solutions of the eigenvalue equation 
for $L$ decaying at plus and minus spatial infinity, the Evans
function vanishes at $\lambda$ if and only if there exists a
solution of the eigenvalue equation decaying at both infinities,
i.e., $\lambda$ is an eigenvalue.
More precisely, zeroes of $D$ agree in both 
location and multiplicity with eigenvalues of $L$ [GJ.1--2].

For traveling waves, there is always an eigenvalue at $\lambda=0$,
corresponding to translational invariance of the underlying PDE,
which in the simplest setting is multiplicity one.
In this situation (which will be the case here), $D(0)=0$,
but $D'(0)\ne 0$.
Moreover, if the evolution equation is well-posed in the sense
that the linearized operator about the wave generates a $C^0$
semigroup, then (by standard resolvent estimates [Pa]) 
$D$ cannot vanish for $\lambda$ real and sufficiently large.
Since $D$, properly constructed respects complex symmetry,
$\bar D(\lambda)=D(\bar \lambda)$, as does the eigenvalue equation
itself, we have in particular that $D$ is real-valued when
restricted to the real axis; thus, $\sgn D(\lambda)$
has a well-defined limit as $\lambda \to +\infty$ along the
real axis, which we will denote as $\sgn D(+\infty)$.
The stability index is then defined as
$$
\Gamma:= \sgn D'(0)D(+\infty).
\eqnlbl{preview}
$$

Evidently, the stability index detects the parity of the number of real roots
on the nonnegative real axis, $\Gamma$ positive corresponding
to even parity, and $\Gamma$ negative to odd parity.
Since complex roots appear in conjugate pairs, this is in fact
the parity of the number of all roots (real and complex) in
the unstable half-plane $\R \lambda >0$, and thus gives
partial information on spectral stability, defined as nonexistence
of eigenvalues in this region.
In particular, $\Gamma >0$ is seen to be a {\it necessary condition
for stability}.
The value of this index comes from the fact that it can be related
to geometric information about the phase portrait of the traveling-wave
ODE, which fact ultimately derives from the
correspondence at $\lambda=0$ between the eigenvalue equation 
and the linearized traveling-wave ODE.

To obtain a concrete result, of course, requires information about
the existence problem, and this can in general be obtained only
in simple situations, e.g. scalar reaction--diffusion equations,
$2\times 2$ parabolic conservation laws, $3\times 3$ conservation
laws with real viscosity, etc., for which the connection problem
is a planar dynamical system, or else [J.1,AGJ,etc.] 
in some singular limit for which the connection problem can be
broken into separately computable ``fast'' and ``slow'' problems.

In this paper we show, in the simple
context of the Majda model, that the methods introduced in
[GZ,BSZ] for the study of stability of viscous shock waves,
may, with slight modifications, be applied also in the study
of stability of detonations to: (i) construct an analytic Evans
function on the set $\R \lambda \ge 0$ for the linearized operator
about the detonation wave, and (ii) in both the strong and
weak detonation cases, compute an expression for $\Gamma$
in terms of quantities associated with the traveling-wave ODE.
In this simple setting, the connection problem is planar, and
we can in fact do more, obtaining a complete evaluation of $\Gamma$;
the result, for both weak and strong detonations, is $\Gamma>0$,
consistent with stability.

This is consistent with prior results of [L,LLT,LYi] for small-amplitude
strong detonations in the small-$q$ limit, 
and of [RV] for arbitrary
amplitude strong detonations in the small-$k$ (ZND) limit.
It is also consistent with results on weak detonations
of the Majda--Rosales model obtained in [Sz] for small-amplitude
waves with intermediate $k$ and in [LY] for arbitrary amplitude waves 
in the large-$k$ limit; however, as far as we know, ours is the
first analytical result on stability of weak detonations for
the Majda model.

These are partial stability results in that they do not rule out
instability;  on the other hand, they are of general applicability,
and are obtained using a relatively
small amount of information about 
the system under study.
In particular, they generalize to the full reactive Navier--Stokes
case (see discussion below), whereas 
analyses of the
integro-differential Majda--Rosales model clearly do not.
Moreover, 
they do have the interesting implication
that transition to instability, if it occurs, must result from a
pair of complex conjugate eigenvalues crossing the imaginary axis,
typically signaling a Poincar\'e--Hopf bifurcation to a time-periodic
solution,
consistent with the experimentally and numerically observed
phenomenon of ``galloping'' detonations 
[MT,FW,MT,AlT,AT,F.1--2, p. 161,BMR,S,Li.6].
It would be very interesting to search numerically for such instabilities
in the high activation-energy limit $k\to \infty$; see, e.g.,
[Br.1--2,BrZ] for an efficient numerical algorithm.

In the course of our development,
we also discuss existence of deflagrations
for a modified ``bump-type'' ignition function, and the appearance
in the ZND limit of degenerate, spatially subalgebraically decaying
families of weak and strong detonations, as described for the
full equations in [GS.1].
Stability of these waves, and their significance in Riemann solutions,
are discussed in Sections 8 and 4, respectively.

\medskip
{\bf Discussion and open problems.}
In [Ly,LyZ],  the methods of this paper are
extended to the full Navier--Stokes
equations of reacting flow, with possibly multi-species reaction.
In this case, the connection problem is high-dimensional, and the
actual evaluation of the geometric stability conditions is carried out
analytically only in the small-$k$ (ZND) limit, for which the 
connection problem has been closely studied in [GS.1--2].
However, the geometric condition itself is valid for arbitrary model
parameters, and could be studied numerically at the same time as
the connection problem; thus, the results again reduce the question
of stability to state of the art existence theory.
Extensions to multi-dimensional stability, in the spirit of viscous
shock calculations of [ZS], are given in [Z.3,JLy.1--2].

As pointed out in [Z.3], Appendix A.3, it is relatively straightforward
using the methods of [ZH,Z.2--4,MZ.1--3] to show that spectral stability
of nondegenerate detonation waves implies 
linearized and nonlinear orbital stability, both for the Majda
model \footnote{For the Majda model with $D=1$, as pointed out in [LLT],
this result may be obtained much more simply by
Sattinger's method of weighted norms [Sat].}
and for the full, reacting Navier--Stokes equations;
this will be a topic of future work.
A very interesting open problem is to carry out a complete spectral
stability analysis in the spirit of the matched asymptotic analysis of the 
Majda model in [RV] for the full reactive Navier--Stokes
equations in the ZND limit, both in one- and multi-dimensions.
The Majda model has the simplifying feature of scalar kinetics;
recent singular perturbation-type techniques developed in the 
shock wave context in [PZ,FreSz] may be helpful in attacking 
the full, system case.
As discussed in Section 8, it would also be interesting to study
further the stability of (generically appearing) degenerate 
weak deflagrations and of degenerate weak and strong detonations
appearing in the ZND limit.

In actual detonations, multidimensional, geometric effects become
important; for example, a converging detonation wave is destabilized,
while diverging, or expansive waves are stabilized, and this is
another important feature to understand; see, e.g. [MR,B,BM,M,Li.7].  
There is also the problem of trying to understand the 
bifurcation to quasi-steady behaviors, such as time-oscillatory
``galloping'', or ``spinning'' detonations in the multidimensional case,
in terms of spectral information given by the Evans function.
The challenge of this problem, similarly as the original stability analysis,
comes from the fact that there is no spectral gap between neutral
point spectrum and essential spectrum of the linearized operator about
the wave, so that tbifurcations are of mathematically nonstandard type;
see the appendix of [BMR] for an interesting related discussion 
in the context of the ZND model.
Useful survey of the latter two topics may be found in [FD,LS,S,Li.6].
%
Finally, we mention the (presumably numerical) problem of 
cataloguing possible Riemann solutions within the class of
viscous profiles in the more general situations discussed
in Section 9
of multi-species reactions or reaction-dependent equation of state;
this appears to be an interesting and physically important
direction for further study.
 
\medskip
{\bf Plan of the paper.} 
In Section 2, we give a brief description of the Majda model,
and the associated Chapman--Jouget diagram.
In Sections 3--4, we discuss existence of profiles and the implications
for solutions of the Riemann problem.
In Sections 5--7 we construct an Evans function for the
linearized operator about the wave,
and carry out the described analysis of the stability index.
In Section 8, we discuss the case of degenerate, spatially
subalgebraically decaying profiles and in
Finally, in Section 9, we describe extensions to multi-species
reactions and reaction-dependent equation of state.

\bigskip
\newsection {\bf Section 2. The Majda model}
\sectionnumber=2 
\theoremnumber=0
\equationnumber=0
\smallskip
\TagsOnLeft
Hereafter, we restrict attention to the Majda model
$$
\cases
u_t + f(u)_x = u_{xx} +kq\varphi(u)z,\\
z_t = -k\varphi(u)z,
\endcases \eqnlbl{1}
$$
or, equivalently,
$$
\cases
(u+qz)_t + f(u)_x =u_{xx},\\
z_t = -k\varphi(u)z,
\endcases \eqnlbl{1alt}
$$
a simplified, scalar version of \eqnref{1abstract}
simulating the dynamics of one-dimensional combustion
within a single characteristic family [M.4].
Here, $u \in \BbbR^1$ is a lumped gas-dynamical variable
combining aspects of specific volume, particle velocity and temperature,
and $z$ corresponds to mass fraction of reactant as before.

Following [M.4], we take
$$
f'(u)>0,\quad f''(u)>0,
\eqnlbl{2.1}
$$
and $q>0$, corresponding to an exothermic reaction.
However, in place of the ``step-type'' function of [M.4],
we take a modified, ``bump-type'' ignition function
$$
\varphi\in C^1,\quad 
\cases
\varphi(u)=
0 &\hbox{for }\quad u\le u_i \hbox{ or } u\ge u^i,\\
\varphi(u)>0
&\hbox{for}\quad u_i< u < u^i.
\endcases \eqnlbl{2.2}
$$
This choice is motivated by the physical parametrization of
temperature with respect to velocity $u$ in the traveling-wave
phase portrait of the ZND model, and
ensures that the traveling-wave equations for \eqnref{1}
agree with the reduced system obtained in [GS.1] by asymptotic
analysis in the ZND limit; see [GS.1], pp. 979--981.
In particular,  it allows for existence of weak deflagration
profiles, as the step-type ignition function does not; see 
Section 3, below.

More precisely, interpreting $u$ as particle velocity we take
$$
\varphi=\psi(T(u)),
\eqnlbl{bump}
$$
where $T(u)$ denotes temperature, 
$\psi$ is a standard, step-type ignition function,
and $\psi\equiv 0$ below some ignition temperature $T_i$ and
positive above, and $T(u)$ is quadratic, concave-down, with
$T(u_i)=T(u^i)=T_i$. 
This agrees qualitatively with the physical dependence of
temperature on velocity along the one-dimensional flow
of the traveling wave ODE for the ZND model.
For example, for an ideal gas with gas constant $\gamma > 1$, 
the physical dependence is
$$
T(u)=-\gamma M^2 u^2 + (\gamma M^2 +1)u,
\eqnlbl{ZND}
$$
where $M=|u_+-s|/|c_+|$, $0<M<1$, 
denotes Mach number of the specific detonation
under consideration, $c_+$ the sound speed at the righthand
state; see [LyZ], eq. (1.17).

\medskip

We seek traveling waves $(u,z)=(\barU(x-st), \Zbar(x-st)$ connecting
states $(u_-,z_-)$ and $(u_+, z_+)$, with 
$$
z_+:= z(+\infty)=1,\quad  z_-:= z(-\infty)=0,
\eqnlbl{2.3}
$$
and
$$
u_i < u_- := u(-\infty) <u^i; \quad u_+:=u(+\infty) \le u_i \hbox{ or }
\ge u^i,
\eqnlbl{2.4}
$$
hence
$$
\varphi_- :=\varphi(u_-) > 0; \quad 
\varphi_+ :=\varphi(u_+)=0:
\eqnlbl{2.5}
$$
that is, combustion waves moving from left to right, leaving
completely burned gas 
in their wake.

\medskip
{\bf The Traveling-wave equation}.  The traveling wave ODE for \eqnref{1alt}
is
$$
\cases
u'=f(u)-f(u_-)-sqz-s(u-u_-),\\
z'=(k/s)\varphi(u) z.
\endcases \eqnlbl{3}
$$

Thus, necessary conditions for existence of a profile are the
modified Rankine--Hugoniot condition 
$$
[f]= s([u]+q),
\eqnlbl{6}
$$
along with
$$
\varphi_+:\varphi(u_+)=0,\quad \hbox{i.e.}\quad u_+ \le u_i \hbox{ or }
\ge u^i,
\eqnlbl{7}
$$
together assuring that $(u_+,z_+)=(u_+,1)$ is a rest state for
\eqnref{3}.

\medskip
{\bf Types of Waves}.  Examining Figure 1, cases A and B, below, 
we find that the possible solutions of \eqnref{6} may be described as follows.

\proclaim{Proposition \thmlbl{1}}  For fixed $u_+$, $s>s_*(u_+)$, there
exist two states $u_->u_+$ for which \eqnref{6} is satisfied.
For $s=s_*$, there exists one solution.  For $s< s_*$, there exist
no solutions $u_->u_+$.  Here, $s_*>f'(u_+)$.
\endproclaim

\proclaim{Proposition \thmlbl{1.1}}  For fixed $u_+$, $s<s^*(u_+)$, there
exist two states $u_-<u_+$ for which \eqnref{6}
is satisfied.
For $s=s^*$, there exists one solution.  For $s> s^*$, there exist
no solutions $u_-<u_+$.  Here, $s^*<f'(u_+)$.
\endproclaim

\midinsert
\captionwidth{2.5truein} 
\hskip .05truein 
\hbox{
\vbox to 2.7truein{\hsize=2.5truein\efig detfig y2.2
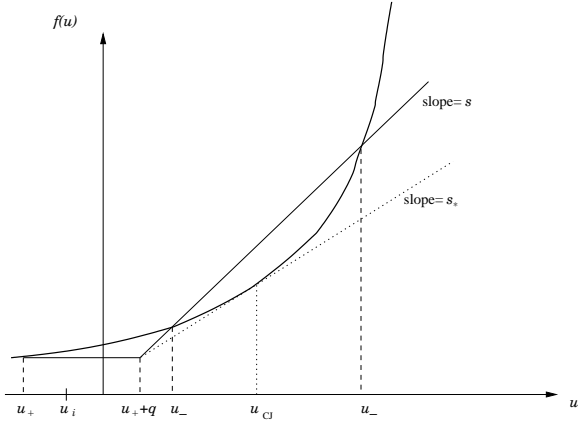
\botcaption{Figure 1a} 
Rest states $u_->u_+$ (detonation case).
\endcaption
}
\hskip .6truein 
\vbox to 2.7truein{\hsize=2.5truein\efig deffig y2.2
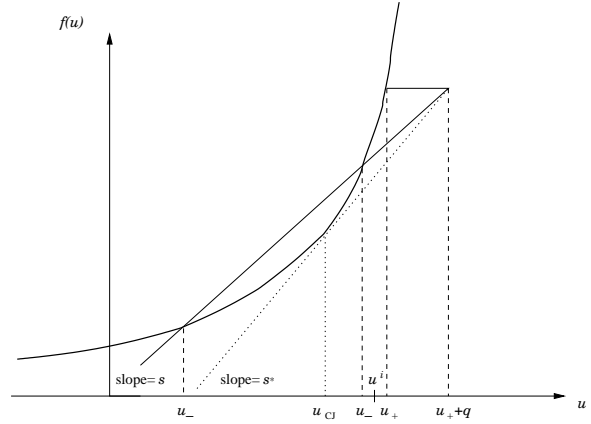
\botcaption{Figure 1b} 
Rest states $u_-<u_+$ (deflagration case).
\endcaption
}
} 
\endinsert

\medskip
Figure 1, known as a {\it Chapman--Jouget diagram},
is the basis for the standard classification of combustion waves.
Denote
$$
a_\pm := f'(u_\pm).
\eqnlbl{8}
$$

{\bf Definition:}  A combustion wave $(\barU,\Zbar)$ is called a 
{\it detonation} if $u_->u_+$.  It is called a {\it strong
detonation} if the Lax characteristic condition holds,
$$
a_->s>a_+.
\eqnlbl{9}
$$
This corresponds to the larger of the two solutions in
Proposition \thmref{1}.  It is called a {\it weak detonation} if
$$
s>a_-, a_+,
\eqnlbl{10}
$$
corresponding to the smaller of the two solutions.  This type is
undercompressive.  The boundary case
$$
a_- = s>a_+
\eqnlbl{10.1}
$$
is 
called a Chapman--Jouget detonation, and has special 
significance in the theory: specifically, in idealized
circumstances, it is the wave expected to be time-asymptotically selected in 
the ``ignition problem'' of initial data consisting of
a large initializing pulse; see [FD,Li.2--4], or discussion at the
end of this section.

\medskip

{\bf Definition:}  A combustion wave $(\barU,\Zbar)$ is called a 
{\it deflagration} if $u_-<u_+$.  It is called a {\it weak deflagration} if
$$
a_+,a_->s,
\eqnlbl{11}
$$
corresponding to the larger of the two solutions, and {\it strong
deflagrations} if
$$
a_+ >s> a_-,
\eqnlbl{12}
$$
corresponding to the smaller.  The former is again undercompressive,
the latter of ``reverse-Lax'' type.
The boundary case
$$
a_+>s=a_-
\eqnlbl{10.1}
$$
is 
called a Chapman--Jouget deflagration.
\medskip


The minimum detonation speed $s_*$ and the maximum deflagration
speed $s^*$ are called Chapman-Jouget (CJ) speeds.
Note that they determine a minimum strength for combustion waves
connecting to $u_+$.  Indeed, there is a band of speeds $(s_*,s^*)$
about the characteristic speed $f'(u_+)$ for which neither detonation 
nor deflagration connections exist, in sharp contrast to the inert-gas, 
shock wave case.
Thus, for fixed $q>0$, 
all combustion waves are ``strong'' in the sense that they differ
appreciably from acoustic signals.

\bigskip
\newsection {\bf Section 3. Existence of profiles.}
\sectionnumber=3 
\theoremnumber=0
\equationnumber=0
\smallskip
\TagsOnLeft

The Chapman--Jouget classification scheme concerns
possible endstates that may be connected by a traveling wave.
We next consider the question of existence of a traveling-wave
profile, reviewing and extending results of [M.4,Sz,Li.3,LY].
\medskip

{\bf Linearized rest points.}  
Linearizing \eqnref{3} about the
critical points $(u_\pm, z_\pm)$, we obtain
$$
\binom {u}{z}' = 
\pmatrix
\alpha _\pm & -sq \\
0 & (k/s) \varphi(u_\pm)
\endpmatrix\ \binom{u}{z},
\eqnlbl{8}
$$
where
$$
\alpha _\pm := a_\pm - s.
\eqnlbl{9new}
$$
(Note:  because of the structure of \eqnref{3}, we cannot without loss of generality set
$s=0$ as in the conservative, shock case.)  
Here, we have used the fact that 
$d\varphi(u_+)=z_-=0$, so that $d\varphi(u_\pm)z_\pm=0$.
Thus, the eigenvalues associated
with the linearized equations are 
$$
\gamma _\pm = \alpha _\pm,\quad (k/s)\varphi(u_\pm).
\eqnlbl{10}
$$
There are two distinct cases: namely, the
generic case that 
$u_+$ differs from ignition temperatures $u_i$ and $u^i$,
i.e., the inequalities in \eqnref{2.5}, \eqnref{7} are strict,
and the degenerate case that
$u_+$ lies precisely at ignition temperature $u_i$ or $u^i$.

In the generic case, referring to \eqnref{2.5}, \eqnref{7}, we find that 
$(u_+, z_+)$ is a {\it saddle--attractor} in the detonation case, 
and a {\it saddle--repellor} in the deflagration case,
with center manifold in both cases lying in 
direction $(qs,\alpha_+)$ tangent to the null-cline $u'=0$
(recall, $\varphi_+=0$), while
$(u_-,z_-)$ is a {\it repellor} in the strong detonation or weak 
deflagration case,  
a {\it  saddle} in the weak detonation or strong deflagration case.
In the exceptional,  Chapman--Jouget case $\alpha_-=0$ that weak and strong
values for $u_-$ coalesce, 
$(u_-,z_-)$ becomes a {\it saddle--repellor},
with center manifold lying in the gas-dynamical direction $(1,0)$;
see \eqnref{8}.
Considering the traveling wave equation as a scalar ODE
$u'=f''(u_-)u^2 + O(|u|^3) -sqz$, 
forced by the solution $z=z_0e^{(k/s)t}$ of
the decoupled exponential growth equation $z'=(k/s)z$,
we find that solutions are generically governed by
the dominant part $u'=f''(u_-)u^2$, growing algebraically as $c/|x|$;
the sole exception is the unique solution lying along the unstable
manifold, $-sqz\sim (k/s)u$, for which both $u$ and $z$ exhibit
exponential growth with rate $e^{(k/s)t}$.

Finally, note that the center manifold (CM) at $(u_+,z_+)$ is a
manifold of {\it rest points}, hence there are no orbits connecting to
$(u_+,z_+)$ along the CM.
That is, for $u<u_i$ or $u>u^i$, ODE \eqnref{3} reduces to
$$
\cases
u'=f(u)-f(u_-)-s(u-u_-)-sqz,\\
z'=0,
\endcases
$$
a {\it scalar} ODE with parameter $z$, for which $u_+$ is
a nondegenerate rest point by the fact that $f'(u_+)\ne s$,
a consequence of Propositions \thmref{1}--\thmref{1.1}:
an attractor in the detonation case (with a unique incoming orbit), 
a repellor in the deflagration case (with no incoming orbit).

\medskip
Collecting the above discussion, we have the preliminary observation:

\proclaim{Lemma \thmlbl{1.1gen}}  In the generic situation
$u_+\ne u_i$, $u_+\ne u^i$, 
detonation profiles, if they exist, are unique up to translation,
decaying exponentially to $(u_+,z_+)$ as $x\to +\infty$.
Moreover, weak and strong detonation profiles
decay exponentially to $(u_-,z_-)$
as $x\to -\infty$,
while Chapman--Jouget detonation profiles generically decay algebraically,
as $c/|x|$.
For fixed $u_+$, $s$, there may exist a weak or a strong detonation,
but not both.
Neither weak nor strong deflagration profiles exist.
\endproclaim

In the degenerate case that $u_+=u_i$
or $u_+=u^i$, the center manifold of rest point $(u_+,z_+)$ 
consists of equilibria only on the ``upper'' side $z>z_+$;
indeed, it is not difficult to see that flow along the CM
is {\it attracting} in either the detonation ($u=u_i$) 
or deflagration ($u=u^i$) case, with rate of approach slower
than any algebraic order. 
For, approximating $z\sim 1$ and $u-u_+\sim (z-z_+)(sq/\alpha_+)$,
we may estimate the order of approach by consideration of the
scalar ODE
$$
(z-z_+)'=(k/s)\varphi(u_+ + (sq/\alpha_+)(z-z_+))\ge 0.
$$
In the detonation case $\alpha_+<0$, for example, this becomes
$\tilde z'=(k/s)\varphi(u_i + (sq/|\alpha_+|)|\tilde z|)$,
where $\tilde z:=z-z_+$.
Noting that $\varphi(u_i+ h)=o(h^k)$ for any algebraic order $k$,
we obtain the result by comparison with ODE $\tilde z'= |\tilde z|^k$
for $\tilde z\le 0$.
In particular, deflagrations now become possible in principle,
while detonations now become possibly nonunique.
More precisely, we have:

\proclaim{Lemma \thmlbl{1.1deg}}  In the degenerate case
$u_+= u_i$ 
strong detonation orbits, if they exist, occur as a one-parameter
family of subalgebraically decaying orbits
bounded on the upper side by a unique exponentially decaying
strong detonation, and on the lower side by a pair of orbits consisting
of a gas-dynamical shock orbit followed by a subalgebraically decaying
weak detonation orbit.
decaying exponentially to $(u_+,z_+)$ as $x\to +\infty$.
Weak detonation profiles, if they exist, are still unique up to
translation, but may be either subalgebraically decaying as $x\to +\infty$,
in which case they bound a one-parameter family of strong detonation orbits
as described above, or exponentially decaying as $x\to +\infty$, in which
case there is no strong detonation profile connecting to $u_+$.
In the degenerate case $u_+=u^i$, there always exists a unique,
weak deflagration profile, decaying subalgebraically as $x\to +\infty$;
strong deflagration profiles, however, do not exist.
The Chapman--Jouget cases are similar, but generically exhibit
algebraic decay $\sim c/|x|$ as $x\to -\infty$, as described in
Lemma \thmref{1.1gen}.
\endproclaim

We defer the proof to the following subsection.
Exponential decay is what is needed to apply the general Evans function
machinery of [GZ].  The ``sonic'' case for shock waves, analogous
to the Chapman--Jouget case has been treated in [H.1--3,HZ] by explicit
calculation, and it appears likely that this approach should generalize
to combustion;  however, we will not treat that issue here, ignoring
for the moment the case of Chapman--Jouget detonations in our 
Evans function analyses.
In the case that $u_+=u_i$ or $u^i$, for 
which
profiles in general approach $(u_+,z_+)$ as $x\to+\infty$ 
at subalgebraic rate, 
it is not possible by current techniques to define the stability index, 
and a different approach must be taken in analyzing stability.

\medskip
{\bf Phase Plane Analysis.}  
A more detailed description of existence properties may be
obtained by examination of the phase plane, as depicted
in Figure 2, below.
Following [LY], consider the nullcline
$$
F(u):= f(u)-f(u_+)-s(u-u_+) -sq(z-z_+)=0
\eqnlbl{F}
$$
on which $u'\equiv 0$.
By convexity of $f$, this intersects the $z=0$ line at
the two points $u_-$ and the $z=1$ line
at the two points $u_+$ described in 
Propositions \thmref{1} and \thmref{1.1};
from here forward, denote these by $u_+< u_-<u^-<u^+$.
By positivity of $sq$, the nullcline $F= 0$ is a graph over $u$,
whence on the physical, invariant region $z\ge 0$,
the region $F\le 0$ above the nullcline is attracting
(recall, $z'\ge 0$ for $z\ge 0$)
so that any flow entering this region for $u>u_i$ remains there
and (since $u'=F<0$) eventually crosses the vertical half-line $u=u_i$, $z>0$.
The value $z=\hat z$ at which it strikes at $u=u_i$ will be the
terminal $z$-value, since $z'\equiv 0$ for $u\le u_i$, and the terminal
$u$-value $\hat u$ will be the unique point $\hat u<u_i$ at which the nullcline
intersects $z=\hat z$.

\midinsert
\captionwidth{2.5truein} 
\hskip 1.5truein 
\hbox{
\vbox to 3.3truein{\hsize=2.5truein\efig phase y2.3
\botcaption{Figure 2} Phase portrait.
\endcaption
}
}
\endinsert

>From these basic considerations, we may conclude already that strong
deflagrations are under no circumstances possible, since these require
a passage from $u_-$ to $u^+$, and these points are separated on $z\ge 0$
by the region $F\le 0$ from which orbits can escape only to the left. 
Weak deflagrations $u^-\to u^+$, if they exist, 
must lie entirely within
the complement $F>0$, hence are monotone increasing in both $z$
and $u$.  It is easily seen that they in fact exist if and only
if $u^+=u^i$ (recall, we have already shown that they do not exist
for $u^+>u^i$).  More generally, there always exists a connection
from $(u^-,0)$ to the rest point $(u^i,z^*)$ determined by the 
intersection of the nullcline $F= 0$ and the vertical
line $u=u^i$.  For, tracing backward along the center manifold 
leading to $(u^i,z^*)$, we see that it is trapped in $F\ge 0$, $z\ge 0$, 
hence must terminate at the corner $(u^-,0)$ of this wedge.

Likewise, weak detonations $u_- \to u_+$, if they exist, are monotone
in both $z$ (increasing) and $u$ (decreasing), lying entirely within
$F\le 0$.  For, the unique orbit originating from saddle $u_-$ into $z\ge 0$
lies along the unstable manifold, pointing into the region $F\le 0$,
from which it can exit only at a rest point along the lefthand edge
of the nullcline $F= 0$.  
Moreover, denoting by $\hat z$ the value of $z$ at which this special 
orbit strikes $u_i$, we find that:
(i) there exists a weak detonation profile if and only if $\hat z=1$,
necessarily unique.
(ii) there exists a unique strong detonation profile if and only if $\hat z< 1$.
(iii) there exists a one-parameter family of degenerate strong 
detonation profiles if and only if $u_+=u_i$, $\hat z=1$, 
and the weak detonation profile approaches $(u_+,1)$ along the center manifold;
if the weak detonation profile approaches along the stable manifold, then
there are no other connections.
(Note that $\hat z\ge 1$ when $u_+=u_i$, since the nullcline
is a lower barrier for the unstable manifold.)
(iv) for $\hat z>1$,  there exist no detonation profiles, 
neither weak nor strong. 

For, in case (iv), this orbit (the unstable manifold from $(u_-,0)$)
serves as a barrier separating $(u_+,1)$ from all orbits in 
$z>0$ originating from $(u^+,0)$.  This is also the case when
there exists a weak detonation profile connecting to $(u_+,1)$
along the stable manifold.  If $\hat z\le 1$ and the unstable manifold
of $(u_-,0)$ does not coincide with the stable manifold of $(u_+,1)$,
then the latter must lie entirely above the former; tracing it backwards,
we find that it must either approach the rest point $(u^-,0)$
within $F\le 0$, corresponding with a monotone strong detonation profile,
or else exit $F\le 0$ along its righthand edge, in which case they must
approach $(u^-,0)$ from within the wedge $F\ge 0$, $z\ge 0$, by the
same argument used to prove existence of weak deflagrations. 
If, along with this nondegenerate strong detonation profile,
there exists also a degenerate weak detonation profile entering
$(u_+,1)$ along a center manifold, i.e., $u_+=u_i$ and $\hat z=1$,
then trapped between these two orbits is a one-parameter family
of degenerate strong detonation profiles also entering along the center
manifold.
This verifies claims (ii)--(iv), while claim (i) is evident.

We have now verified the claims of Lemma \thmref{1.1deg};
more, we have characterized the existence of weak and strong
detonation profiles in terms of the height 
$\hat z$ at which the unstable
manifold of $(u_-,0)$ strikes $u_i$
(and in degenerate cases also the direction): 
equivalently, in terms of the {\it Melnikov separation function}
$$
d(u_+,s,k,q):= \hat z-1.
\eqnlbl{melnikov}
$$
The following monotonicity properties will be helpful
in completing our description of existence.

\proclaim{Proposition \thmlbl{mon}.}
The function $d(\cdot)$, 
more generally, the height of the unstable manifold of $(u_-,0)$ 
at any value $u_i\le u < u_-$,
is monotone nondecreasing in $u_+$, $k$, and, for $\hat z\le 1$, in $q$,
and is nonincreasing in $s$.
Moreoever, this monotonicity is strict except 
in the special case that  $(u_i,\hat z)$ is a rest point,
i.e., lies on the nullcline $F=0$.
\endproclaim

{\bf Proof.}
To establish monotonicity, we show in each case that the vector field
$(F, G)$, $G:= (k/s)\varphi(u)z$, determining the flow of the traveling wave ODE,
within the region $F\le 0$, $z\ge 0$, $u_i<u < u_-$
rotates clockwise {\it strictly} monotonically with respect to the 
parameter $\beta$ under consideration, or, equivalently,
$$
\delta_\beta:=
\det 
\pmatrix F & \partial F/\partial \beta \\
G & \partial G/\partial \beta \\
\endpmatrix \lessgtr 0,
\eqnlbl{mon}
$$
with $<0$ corresponding to monotone increase, and $>0$ to monotone decrease.
For example, $\delta_k$ is simply $FG/k<0$,
$\delta_q$ is $Gs(z-z_+)<0$, and $\delta_{u_+}= G \alpha_+<0$,
by the fact that  $\alpha_+< 0$ (weak detonation).
Finally, 
$$
\aligned
\delta_s:=
\det 
\pmatrix F & -q(z-z_+)-(u-u_+)\\
G & -(k/s^2)\varphi(u)z\\
\endpmatrix 
&=
\det 
\pmatrix F & (1/s)(F- (f(u)-f(u_+))\\
G & -(1/s)G\\
\endpmatrix \\
&=
(-G/s)(2F- (f(u)-f(u_+)) \\
& > 0.
\endaligned
\eqnlbl{deltas}
$$

This gives strict comparison principles on $u_i < u<u_-$,
not only for the unstable manifolds, but for any orbits
lying in $F\le 0$, $z\ge 0$ ($1\ge z\ge 0$ in the case of $q$),
$u_i<u<u^i$.
Noting that $\partial u_-/\partial s <0$, 
$\partial u_-/\partial q$,
$\partial u_-/\partial u_+ > 0$,  and
$\partial u_-/\partial k=0$, we thus obtain the claimed
strict monotonicity on $u_i<u<u_-$.
It follows that there is an orbit lying strictly between
two stable manifolds under consideration for $u_i<u<u_-$,
and this implies strict monotonicity at $u=u_i$ unless this
third orbit collides with the stable manifold at $u_i$, which
can only happen if $(u_i,z_*)$ is a rest point.
\myqed

(Note: the above argument applies also in the Chapman--Jouget case.)

\medskip
{\bf Remark \thmlbl{mderiv}.}
In the nondegenerate case $u_+\ne u_i$,
we may express the derivatives of $d$ at a
profile $d=0$
in the standard way (see, e.g., [GH,HK]) as {\it Melnikov integrals} 
$$
\partial d/\partial \beta=
-\int^\infty_{-\infty}
e^{-\int_0^x \Tr(\partial(F,G)/\partial (u,z))(\bar u, \bar z)(y)) dy}
\det
\pmatrix
F &  \partial F/\partial \beta \\
G &  \partial G/\partial \beta \\
\endpmatrix(\bar u, \bar z)(x) dx,
\eqnlbl{mel}
$$
to find by the same calculation that the derivatives
$\partial d/\partial \beta$ have strict signs as well:
$\partial d/\partial u_+$,
$\partial d/\partial k$, and $\partial d/\partial q>0$,
while $\partial d/\partial s<0$.
In particular, we have
$$
\aligned
\partial d/\partial s&=
\int^\infty_{-\infty}
e^{-\int_0^x(\alpha+ k\varphi)(\bar u(y)) dy}
\det
\pmatrix
\barU_x &(\barU-u_+)+ q\bar (z-z_+)\\
\Zbar_x &(1/s)\Zbar_x
\endpmatrix(x) dx\\
&=
\int^\infty_{-\infty}
e^{-\int_0^x(\alpha+ k\varphi)(\bar u(y)) dy}
(\bar z_x/s)
(2\bar u'-(f(\bar u)-f(u_+)))(x)
\, dx
<0.\\
\endaligned
\eqnlbl{33.2}
$$
\medskip
{\bf Remark \thmlbl{u_-}.}
Alternatively, following [M.4], we may fix the lefthand state state $u_-$ 
and consider $d$ as a function of $(k,q,u_-,s)$.
It is a straightforward exercise to verify that $d$ then becomes monotone
nondecreasing in $u_-$ and  $k$, but {\it nonincreasing} in $q$
(now for all values of $\hat z$).
Monotonicity in $s$ is lost.
\medskip

%

Our main conclusions regarding existence are summarized in the 
following two propositions.

\proclaim{Proposition \thmlbl{5.1}}  For each $u_+$, there holds
one of two possibilities: either (i) for each speed $s$ 
greater than or equal to the minimal, Chapman--Jouget detonation speed 
$s_*$ for which detonations can occur, there exists a single exponentially
decaying strong detonation profile,
or (ii) there exists a threshold speed $\hat s\ge s_*$ at which there
exists a (necessarily) unique exponentially decaying weak detonation
profile, above which there exists a single exponentially
decaying strong detonation profile, and below which there exists
neither weak nor strong profile.
Where it is defined, $s_*$ is monotone increasing in $u_+$.
If $u_+\ne u_i$, then strong detonation profiles, when they exist,
are the only detonation profiles that occur. 
If $u_+=u_i$, on the other hand, then whenever there exists an exponentially
decaying strong detonation,  it bounds from above a degenerate family of
strong detonation profiles, bounded below by a degenerate weak detonation
profile and a gas-dynamical shock profile,
as described in Lemma \thmref{1.1deg}.
Weak detonations are always monotone in $u$ (decreasing) and $z$
(increasing).

Weak deflagrations exist if and only if $u^+=u^i$, and are 
monotone increasing in both $u$ and $z$;
moreover, they are always degenerate, decaying sub-algebraically
as $x\to +\infty$.
Strong deflagrations do not occur.
\endproclaim

\proclaim{Proposition \thmlbl{ZND}} For fixed model parameters
$f$, $\varphi$, $q$, there exists $k_0>0$ such that, in
Proposition \thmref{5.1},
only case (i) occurs for $0<k<k_0$, while, for $k\ge k_0$,
case (ii) always occurs for some choice of $(u_+,s)$.
\endproclaim

\medskip
{\bf Proof of Proposition \thmref{5.1}.}
It is straightforward to see, in the infinite-speed limit
$u_-\to u_++q$, that $\hat z <1$ if $u_+\ne u_i$, and
$\hat z=1$ if $u_+=u_i$, with the weak detonation profile 
entering $u_+=u_i$ along the center manifold.
If $u_++q<u_i$ is out of range, then we consider instead the
maximum speed limit $u_- \to u_i$, observing that $\hat z\to 0$
is again less than $1$.
As speed $s$ is decreased, $\hat z$ increases monotonically,
so that either $\hat z\le 1$ for all $s$, with the unstable
manifold of $u_-$ lying always below the stable manifold of
$u_+$, or else there is a unique transition point $s_*$ as
described in the theorem at which the two manifolds coincide,
below which $\hat z>1$ and there are no detonation connections.
That $s_*$ is monotone in $u_+$ follows from the fact that
$\hat z$ decreases with $s$ but increases with $u_+$.
The remaining conclusions about detonations, and the conclusions
about weak and strong deflagrations, follow from the discussion
above Proposition \thmref{mon}.
\myqed

\medskip
{\bf Proof of Proposition \thmref{ZND}.}
As $k\to 0$, the unstable manifold of $(u_-,0)$ approaches the nullcline
$F=0$, as can be seen either by direct calculation, or using
singular perturbation techniques as in [GS.1].
Likewise, the stable manifold of rest point $(u_i,z_*)$
(recall: the intersection of the nullcline and $u=u_i$) 
approaches a horizontal line.
Thus,  the unstable manifold of $u_-$ is trapped between
the nullcline and the stable manifold of $(u_i,z_*)$,
hence must approach $(u_i,z_*)$ along the center manifold
tangent to the nullcline, decaying subalgebraically as
$x\to +\infty$; see the discussion of case $u_+=u_i$ for
details.  
(In particular, it connects to $(u_+,1)$ if and
only if $u_+=u_i$.)
>From these observations, and $\hat z=z_*\le 1$, we find that
we are in case (i).
Moreover, this argument is uniform in model parameters, yielding a
global result for $k$ less than some threshold $k_0$ as described.
That case (ii) occurs for all $k\ge k_0$ follows by monotonicity
of $\hat z$ with respect to $k$.
\myqed


%

%

\bigskip
\newsection {\bf Section 4. 
Riemann solutions and the CJ shift.}
\sectionnumber=4
\theoremnumber=0
\equationnumber=0
\smallskip
\TagsOnLeft

We conclude our discussion of existence by cataloging
solutions of {\it Riemann 
problems} involving waves with viscous profiles.
In a Riemann problem, we
prescribe data $U_L=(u_L,z_L)$ for $x\le 0$ and $U_R=(u_R,z_R)$ for $x> 0$,
and
seek a sequence of waves $(U_L,U_1)$, $(U_1,U_2)$, $\dots$,
$(U_m,U_R)$ with increasing speeds,
progressing from left state $U_L$ to right state $U_R$,
consisting of gas-dynamical shock or rarefaction waves,
weak or strong detonations, or weak deflagrations, each
possessing a viscous profile, 
for which the endstates $U_j=(u_j,z_j)$, $j=1,\dots,m, R$ are
individually valid time-asymptotic states, i.e., $z_j=0$ if $u_i<u_j<u^i$.
As discussed, e.g., in [AMPZ], such solutions represent 
possible time-asymptotic states for solutions of \eqnref{1}, \eqnref{1alt},
as obtained by formal, matched asymptotic expansion with the usual,
hyperbolic scaling $x,t \to \varepsilon x,\varepsilon t$.

By the one-sided nature of the Majda model, it is easily seen that
the only interesting situation is the one $z_L=0$, $z_R=1$ described
in the introduction.
For, the reverse situation $z_L=1$, $z_R=0$ clearly admits no solution,
there being no waves along which $z$ increases, nor, for the same
reason, does the situation $z_L=z_R=1$ with $u_L$ and $u_R$
lying in different components of $(-\infty,u_i] \cup [u^i,+\infty)$.
In these cases, the asymptotic behavior has a different, diffusive
scaling not captured by the Riemann solution.  
For $z_L=z_R=0$ on the other hand, or $z_L=z_R=1$ and
$u_L$ and $u_R$ lying in the same component of 
$(-\infty,u_i] \cup [u^i,+\infty)$, the solution is just the
gas-dynamical solution for the conservation law $u_t+f(u)_x=u_{xx}$
with the value of $z$ held fixed.
We therefore fix $z_L=1$, $z_R=0$ for the remainder of our discussion.

%

\medskip
{\it Case I. (deflagration: $u_L< u^i\le u_r$).}
Denote by $u_{CJ}$ the Chapman--Jouget deflagration speed
associated with $u^+=u_R$.  Then, there are two subcases.
{\it Ia.} For $u_{CJ}< u_L < u^i$, the solution consists of a weak
deflagration from $u_L$ to $u^i$, followed by a (possibly zero strength)
fluid-dynamical rarefaction with $z\equiv 1$ from $u^i$ to $u_R$.
{\it Ib.} For $u_L< u_{CJ}$, the solution consists of a
(possibly zero strength) fluid-dynamical rarefaction from $u_L$ to $u_{CJ}$,
followed by a Chapman--Jouget deflagration from $u_{CJ}$ to $u^i$,
followed by a fluid-dynamical
rarefaction with $z\equiv 1$ from $u^i$ to $u_R$.
Note that weak deflagrations may be followed by
gas-dynamical waves, whereas the CJ deflagration may be both
followed by and preceded by another wave, so long as the preceding
wave is a rarefaction.

\medskip
{\it Case II. (standard detonation: $u_L< u^i$, $u_r< u_i$).}
In this case, the solution structure depends on whether
we are in case (i) or case (ii) as described in Proposition \thmref{5.1}.
Denote by $u_{CJ}$ the Chapman--Jouget detonation speed
associated with $u^+=u_R$.  
If we are in case (i), then there are two subcases.
{\it IIa.} For $u_{CJ}< u_L< u^i$, the solution consists simply of
a strong detonation from $u_L$ to $u_R$.
{\it IIb.} For $ u_L\le u_{CJ}$, the solution consists of
a (possibly zero strength) gas-dynamical rarefaction from
$u_L$ to $u_{CJ}$, followed by a Chapman--Jouget detonation from
$u_{CJ}$ to $u_R$.

If we are in case (ii), there are again two subcases, but with the
special role of the CJ detonation, which no longer has a profile,
now played by the unique weak detonation possessing a (nondegenerate)
profile.  Denote by
$u_* \le u_{CJ}$ the special left state $u_-$ that is connected to $u_R$
by a nondegenerate weak detonation profile with speed $s_*$,
and $u^* \ge u_{CJ}$ the corresponding value of $u^-$: i.e.,
the unique state that is connected to $u_*$ by a gas-dynamical
shock with the same speed $s_*$.  
{\it IIa'.}  For  $u^*<u_L < u^i$, the solution consists simply of
a strong detonation from $u_L$ to $u_R$.
{\it IIb'.}
For  $u_L\le u^*$, the solution consists of a (possibly zero strength)
gas-dynamical shock or rarefaction, according as $u_L$ is
$\ge$ or $\le u_*$, from $u_L$ to $u^*$, followed by a weak detonation
from $u_*$ to $u_R$.

{\it Case III. (limiting detonation: $u_L< u^i$, $u_r= u^i$).}
In this degenerate case, there exist, along with each nondegenerate
strong detonation profile from $u^-$ to $u_i$, a degenerate weak detonation 
profile from $u_-$ to $u_i$.
This leads to nonuniqueness of Riemann solutions for $u_L$ on the
interval between $u_{CJ}$ or $u^*$ (in case (i) or (ii), 
respectively) and $u^i$.  Namely, along with the solutions described in case II,
we also have solutions obtained by substituting for any strong detonation
with speed $s> s_*$ from $u^-$ to $u_+$
a gas-dynamical shock with speed $\le s$ from $u^-$ to $\hat u\le u_-$, 
followed by a degenerate weak detonation with speed $\hat s \ge s> s_*$ from
$\hat u$ to $u_+$ (recall: both strong detonation and degenerate weak
detonation profiles exist for any speed $>s_*$).
There is a milder nonuniqueness at the level of profiles, since any
nondegenerate strong detonation profile may be replaced by a degenerate
one with the same endstates.

{\it Case IV. ($u_L\ge u^i$, $u_R> u_i$).}
It is easily deduced that this problem has no solution.
For, in order that $z$ increase from $0$ to $1$, there must
be an intermediate state within $(u_i,u^i)$, and this can
only be reached from $u_L$ by a single gas-dynamical shock.
The only wave that can follow a shock is a weak detonation,
and no wave can follow a detonation.  But, $u_R$, as the right endstate
of a weak detonation must then satisfy $u_R\le u_i$, a contradiction.

{\it Case V. ($u_L\ge u^i$, $u_R<u_i$).}
By the discussion of the previous case, 
the only possible solution is a gas-dynamical
shock from $u_L$ to $u_1$, followed by a weak detonation from
$u_-=u_1$ to $u_+=u_R$ having the special property
$u_+\le u_i < u_- < u^i\le  u_L \le u^-$.
This never occurs in detonation case (i), but may occur in case (ii)
within certain parameter range.



\medskip

Note that the solvable cases I--III contain the basic {\it propagation problem} 
described in the introduction, of a combustion wave moving from a burned,
superignition state at $-\infty$ 
into an unburned, subignition state at $+\infty$.
Likewise, cases II-III and IV contain the {\it ignition problem}
$U_L=(u_L, 0)$, $U_R=(u_R,1)$, where $u_L, u_R \le u_i$ or 
$u_L, u_R\ge u^i$: a one-sided version 
obtained by left-right symmetry 
of the ignition problem $U_L=(u_L,1)=U_R$ 
for the full, reacting Navier--Stokes equations, 
corresponding to a large, pulse-type
excitation, initiating combustion,
of a quiescent background state.
We do not have a useful interpretation of case V.
%
In cases I--II, the Riemann solution is uniquely determined within the
class of viscous profiles, while in case III it is not.  If we restrict
to the class of nondegenerate (exponentially decaying) detonation profiles,
then the Riemann solution is uniquely determined also in case III.
A special role in the solution structure is played by the detonation 
profile with minimum speed:  the CJ detonation, in case (i),
the weak detonation profile with speed $s_*$, in case (ii).

Unique solvability in the class of exponentially decaying profiles
suggests that nondegenerate detonation profiles, both weak and strong, 
at least in usual circumstances are stable.
Change in stability would presumably signal bifurcation to more complicated
time-asymptotic behavior, perhaps involving time-oscillatory
``galloping'' solutions.

\medskip
{\bf Remark \thmlbl{shift}. (CJ shift)}
In both propagation and ignition Riemann problems considered above,
appearance of a weak detonation profile, case (ii), is associated
with a shift in the solution structure from CJ to weak detonation,
in contrast with the behaviour predicted by the ZND model.
Such a shift is indeed observed in experiments, for appropriate
parameter regimes [FD].
Proposition \thmref{ZND} states that no such shift
occurs for $k<k_0$: that is, it validates the ZND picture for
$k$ merely small, and not only in the $k\to 0$ limit.
This important observation was made first in [GS.1], in the
larger context of the full reacting Navier--Stokes equations,
and answers in the negative a conjecture of Majda [M.4].
\footnote{In the Remark below Theorem 1 of the reference, the
conjecture that $q_0^{CR}>\hat q$ for any $k_0>0$: in our notation,
that $\hat z>1$ for $u_+=u_i$, for any $k>0$,
which would imply the existence of nondegenerate weak detonation
profiles for $u_+$ less than but sufficiently close to $u_i$. }

\medskip
{\bf Remark \thmlbl{linwp}.}
As should be apparent from the above discussion, it is monotonicity
of $d$ with respect to $s$ that leads to unique solvability of
the Riemann problem in the presence of weak detonation profiles 
(detonation case (ii)).
Likewise, a local analysis in the spirit of
[ScSh,ZPM,Fre.1,GZ,ZS] shows that
$\partial d/\partial s\ne 0$ is equivalent to
linearized well-posedness of the Riemann problem about weak
detonation data.

\bigskip
\newsection {\bf Section 5. Construction of the Evans function.}
\sectionnumber=5
\theoremnumber=0
\equationnumber=0
\smallskip
\TagsOnLeft

We now turn to the question of stability, 
beginning by a careful construction of the Evans function.
Consider a nondegenerate, i.e., spatially exponentially decaying
traveling wave profile $(\barU(x-st)$ of any type,
satisfying the nonsonicity assumption $a_\pm \ne s$.
The linearized equations of
\eqnref{1}, \eqnref{1alt}
about $\barU(x-st)$, in moving coordinates $\tildeX=x-st$, are
$$
\cases
u_t 
-q( k\varphi '(\barU)u \Zbar - k \varphi(\barU)z)
+ (\alpha  u)_x  = u_{xx}\\
z_t - sz_x = 
-\ k\varphi '(\barU)u \Zbar -
 k \varphi(\barU)z,
\endcases\eqnlbl{13}
$$
$$
\cases
u_t + qz_t + (\alpha  u)_x - sqz_x = u_{xx}\\
z_t - sz_x = 
-\ k\varphi '(\barU)u \Zbar - k \varphi(\barU)z,
\endcases\eqnlbl{13alt}
$$
where
$$
\alpha := a(\barU)-s,\quad a(u) = df(u).
\eqnlbl{14}
$$
The associated eigenvalue equation is thus
$$
\cases
u'' = (\alpha u)'+\lambda u 
-q(k \varphi '(\barU)u \Zbar+k\varphi (\barU) z )
,\\
z' = (1/s) 
(k \varphi '(\barU)u \Zbar+k\varphi (\barU) z + \lambda z),
\endcases \eqnlbl{15alt}
$$
or, alternatively,
$$
\cases
u'' = (\alpha u)'+\lambda u +q(\lambda z-sz'),\\
z' = (1/s) (k \varphi '(\barU)u \Zbar+k\varphi (\barU) z + \lambda z),
\endcases \eqnlbl{15}
$$
and the limiting systems at $x=\pm \infty$ are
$$
\cases
u'' = \alpha _\pm u' +\lambda u -qk\varphi_\pm z,\\
z' = (1/s) (k \varphi_\pm z + \lambda z).
\endcases \eqnlbl{16}_\pm
$$
(Recall: $\varphi_- >0$, $z_-=d\varphi_+=\varphi_+=0$.)
Here, $\alpha_\pm \ne 0$, by the nonsonicity assumption $a_\pm \ne s$.

Seeking solutions $u=e^{\mu x}U,\ z=e^{\mu x}Z$ of \eqnref{16}$_\pm$,
we obtain the characteristic equation
$$
\pmatrix
\mu^2-\mu \alpha _\pm - \lambda  & -qk\varphi_\pm\\
0 & \mu-(k/s)\varphi_\pm - \lambda/s 
\endpmatrix \quad 
\binom{U}{Z}= \binom{0}{0}.
\eqnlbl{17}_\pm
$$
This is readily solved using its block triangular form,
to yield growth rates (eigenvalues)
$$
\mu=\frac{\alpha _\pm  \pm \sqrt{ \alpha _\pm^2+4\lambda }}{2},
\quad 
\underbrace{(k/s)\varphi_\pm}_{\ge 0} +\lambda/s, 
\eqnlbl{18}_\pm
$$
associated, for all except finitely many points $\lambda$ where
eigenvalues may coincide, with respective normal modes of form
$$
\pmatrix u \\u ' \\z\endpmatrix=
e^{\mu x}\pmatrix U\\ \mu U \\Z\endpmatrix=
e^{\mu_\pm x} \pmatrix 1\\\mu_\pm  \\ 0
\endpmatrix,
\quad 
e^{\mu x}
\pmatrix  * \\*\\ 1
\endpmatrix.
\eqnlbl{18.5}_\pm
$$

Inspection of formulae \eqnref{18}$_\pm$--\eqnref{18.5}$_\pm$ yields

\proclaim{Lemma \thmlbl{2}}  On $\{\R \lambda >0\}$, there are {\it
two} unstable and {\it one} stable modes of \eqnref{17}$_\pm$ at either
of $\pm \infty$,
$$
\mu^\pm_1 < 0 < \mu^\pm_2, \mu^\pm_3;
$$
moreover, there exist
analytic choices of solutions of form
$$
\pmatrix u\\u' \\z\endpmatrix(\lambda)=
\pmatrix u\\ \mu_\pm u  \\ 0
\endpmatrix,
\quad 
\pmatrix  * \\*\\ z
\endpmatrix
$$
spanning the associated stable/unstable manifolds.
Moreover, each of these analytic functions has the
complex symmetry $\bar f(z)=f(\bar z)$.
\endproclaim

{\bf Proof.} The first claim is immediate from the formulae;
the second then follows as described in [GZ,Z.3]
by a standard lemma of Kato [K], which asserts that analytic
subspaces on a simply connected domain possess an analytic
choice of bases.
The third statement follows exactly as in [GZ,Z.3] from the observation
that the construction of Kato preserves complex symmetry.
\myqed

Further, we have

\proclaim{Lemma \thmlbl{3}}  The functions $\mu^\pm_j, u^\pm_j$ and 
$z^\pm_j$ can be analytically extended onto a region 
$$
\Omega :=\{\R \lambda \ge -\theta_1-\theta_2|Im \lambda |^2\},
\quad \theta_1, \theta_2 > 0.
\eqnlbl{19}
$$
Moreover, $\mu_j$ are analytic on a neighborhood of $\lambda=0$, 
with associated (distinct) analytic eigenvectors.
\endproclaim

{\bf Proof.} A closer look reveals that Lemma \thmref{2} holds true
on $\Omega \setminus B(0,r)$ for arbitrary $r>0$ and
$\theta_j(r)>0$ sufficiently small.
Thus, it is sufficient to show that there exists an analytic
extension on $B(0,r)$: in particular, to verify the second claim.
At $x=-\infty$, $\varphi_->0$, and so, at $\lambda=0$, 
the reactive root is $\mu= (k/s)+\lambda/s>0$, 
while the kinematic roots are $0$ and $\alpha_+\ne 0$;
thus, there is a spectral gap between the single zero mode
and strong stable or unstable modes, and so each of these may be
analytically continued through the origin.
At $x=+\infty$ on the other hand, 
$\varphi_+=0$, and \eqnref{17}$_\pm$ becomes block
diagonal, decoupling into kinematic and reactive blocks.  Within
blocks, the roots $\mu$ are distinct at $\lambda=0$, and the result follows
as before.
\myqed

{\bf Remark.} For the purposes of this paper, it is sufficient 
to establish Lemma \thmref{3} on the smaller set 
$\Omega= \{\lambda: \, \R \lambda \ge 0\}$. 
\medskip

Applying the Gap lemma of [GZ,KS], \footnote{Here, we make use
of the assumed spatial decay of the wave, a hypothesis of
the Gap lemma.}
we obtain, finally,

\proclaim{Corollary \thmlbl{4}}  There exist solutions $u^\pm_j
(\lambda ,x),\, z^\pm_j(\lambda ,x)$ of eigenvalue equation
\eqnref{15}, analytic on $\Omega$, \eqnref{19}.  Moreover, on $\{\R
\lambda >0\},\, (u^+_1, z^+_1)$  spans the stable manifold at
$+\infty$, and $\{(u^-_2, z^-_2),\, (u^-_3,z^-_3)\}$ the unstable
manifold at $-\infty$.
\endproclaim

{\bf Definition:}  Following the development of
[GZ] we define the {\it Evans function} as 
$$
D(\lambda ):=\det
\pmatrix
u^+_1 &u^-_2 & u^-_3\\
{u^+_1}' &{u^-_2}' &{u^-_3}'\\
z^+_1 & z^-_2 & z^-_3
\endpmatrix_{|_{x=0}} .
\eqnlbl{20}
$$

\medskip
\proclaim{Proposition \thmlbl{evansprops}}
The Evans function $D$ is analytic on $\Omega$; on the
subdomain $\R \lambda>0$, its zeroes
correspond precisely to eigenvalues of the linearized operator $L$. 
Moreover, it possesses complex symmetry
$D(\bar z)=\bar D(z)$; in particular, $D(\lambda)$ is real for
real $\lambda$.
\endproclaim

{\bf Proof.}
Analyticity and complex symmetry are inherited from the properties of
the columns of the matrix from which $D$ is evaluated.
Recalling that, on $\R \lambda >0$, the columns of that matrix
span the decaying manifolds of the eigenvalue equation at $\pm \infty$,
we find, further, that vanishing of $D$ is equivalent to nontrivial
intersection of these manifolds, i.e., existence of a solution decaying
at both spatial infinities.  Thus, zeroes of $D$
correspond in location with eigenvalues of $L$.  That they correspond
also in multiplicity follows from a more detailed calculation of
Gardner and Jones [GJ.1--2]; see also Section 6 of [ZH].
\myqed

\proclaim{Corollary \thmlbl{basicstab}}
A necessary condition for linearized stability is that 
the number of zeroes of $D$ in the unstable half-plane
$\{\lambda: \, \R \lambda >0\}$ be zero; in particular,
that it be even.
\endproclaim
\bigskip
\newsection {\bf Section 6. Stability analysis for strong detonations.}
\sectionnumber=6
\theoremnumber=0
\equationnumber=0
\smallskip
\TagsOnLeft

Specializing now to the strong detonation case, $\alpha_->0>\alpha_+$,
we complete our stability analysis by a computation of
the stability index $\Gamma:=\sgn D'(0)D(+\infty)$.
\medskip

{\bf Analysis at $\lambda =0$}.  At $\lambda =0$, \eqnref{15} reduces to the
linearized traveling wave ODE
$$
\cases
u'' = (\alpha u)' -sqz',\\
z' = (1/s) (k \varphi ' \Zbar u + k\varphi(\barU)z).
\endcases \eqnlbl{21}
$$
Take without loss of generality
$$
(u^+_1,z^+_1) = (u^-_3,z^-_3) = (\barU_x, \Zbar_x),
\eqnlbl{22}
$$
both exponentially decaying modes.  The solution $(u^-_2, z^-_2)$
is likewise asymptotically decaying at $-\infty$,
$$
u^-_2(-\infty)=z^-_2(-\infty)=0,
\eqnlbl{23}
$$
since, for $\alpha_->0$, the unstable modes at $-\infty$
remain strictly unstable as $\lambda \to 0$.

Integrating the first equation of \eqnref{21} from $\pm \infty$ to
$x$, we thus obtain the homogeneous first-order equation
$$
{u}'-\alpha u+sqz = 0,
\eqnlbl{25}_\pm
$$
satisfied by each of $(u_1^+,z_1^+)$, $(u_2^-,z_2^-)$, $(u_3^-,z_3^-)$.

\proclaim{Proposition \thmlbl{4}}  $D(0)=0$, while
$$
D'(0)=\gamma \Delta
$$
where
$$
\gamma =\lim_{M\to +\infty}
\det \pmatrix
u_2^- & \bar u_x \\
z_2^- & \bar z_x 
\endpmatrix
(-M)
e^{\int^0_{-M}(\alpha
(x)+(k/s)\varphi\barU(x))}.
\eqnlbl{25.1}
$$
and 
$$
\Delta=[u]+q.
\eqnlbl{deltadef}
$$
\endproclaim

{\bf Proof}.  $D(0)=0$ is immediate from \eqnref{22}.  Using the
Leibnitz rule, we obtain as usual 
$$
D'(0)=\det
\pmatrix
\barU_x & u^-_2 & y\\
\barU_{xx} & {u^-_2}' &y'\\
\Zbar_x & z^-_2 & \sigma
\endpmatrix , \eqnlbl{26}
$$
where $y:=y^- -y^+,\ \sigma := \sigma^- - \sigma^+,$ and
$(y^+,\sigma^+):= (\partial/\partial \lambda)  (u^+_1,z^+_1),
\ (y^-,\sigma^-):= (\partial/\partial \lambda) (u^-_3,z^-_3)$, satisfy
the variational equations
(obtained by differentiating \eqnref{15})
$$
\cases
{y^\pm}'' = (\alpha y^\pm )'-{sq\sigma^\pm}' + \barU_x + q\bar z_x,\\
{\sigma^\pm}' = (1/s) (k \varphi'(\barU)y^\pm \Zbar + k
\varphi(\barU)\sigma^\pm +\Zbar _x) ,
\endcases
\eqnlbl{27}
$$
with
$$
(y^\pm,\sigma^\pm)(\pm \infty)=0.
\eqnlbl{28}
$$
Thus
$$
{y^\pm}'-\alpha y^\pm + sq\sigma^\pm = \barU-u_\pm + q(\bar z - z_\pm),
\eqnlbl{28}
$$
and
$$
y-\alpha y+sq\sigma=[u]+ q,
\eqnlbl{29}
$$
where $[u]:=(u_+-u_-)$.  Using \eqnref{25}$_\pm$, \eqnref{28} to partially
eliminate the 2nd row in \eqnref{26}, we obtain
$$
\aligned
D'(0)
&=\det
\pmatrix
\barU_x & u^-_2 & y\\
0 &0 &[u]+q\\
\Zbar_x & z^-_2 & \sigma
\endpmatrix_{|_{x=0}}\\
&= - \det
\pmatrix
\barU_x & u^-_2\\
\Zbar_x & z^-_2
\endpmatrix_{|_{x=0}}\ ([u]+q)\\
&=\gamma \Delta,
\endaligned \eqnlbl{30}
$$
where $\gamma := \det
\pmatrix
u^-_2 & \barU_x \\
z^-_2& \Zbar_x
\endpmatrix |_{x=0}$ is a Wronskian for the linearization about $(\barU, \Zbar)$ 
of the traveling wave equation \eqnref{3}.  
The evaluation
\eqnref{25.1} then follows by Abel's formula.\myqed

\medskip
{\bf Analysis as $\lambda \to +\infty$}.  
Similarly as in [GZ], we next establish

\proclaim{Proposition \thmlbl{5}}  As $\lambda \to \infty$ along the
real axis,
$$
\sgn D(\lambda )\to  \sgn u_1^+(+\infty)
\det
\pmatrix
u^-_2 & u_3^-\\ 
z^-_2& z_3^-
\endpmatrix (-\infty)
\eqnlbl{largelambda}
$$
\endproclaim

{\bf Proof.}
Following [GZ], consider the frozen-coefficient version
$$
\cases
u'' = \alpha(\bar u) u'+ \alpha(\bar u)_x u+\lambda u 
-qk\varphi z,\\
z' = (1/s) (k \varphi '(\barU)u \Zbar+k\varphi (\barU) z + \lambda z),
\endcases \eqnlbl{frozen}
$$
of \eqnref{15} obtained by evaluating $\bar u$ at a fixed point $x_0$,
and as usual seek solutions $u=e^{\mu x}U$, $z=e^{\mu x}Z$: that is,
solutions of characteristic equation
$$
\pmatrix
\mu^2-\mu \alpha -\alpha_x  - \lambda  & -qk\varphi\\
-(k/s) \varphi'\bar z 
 & \mu-(k/s)\varphi - \lambda/s 
\endpmatrix \quad 
\binom{U}{Z}= \binom{0}{0}.
\eqnlbl{frozenchar}
$$

It is a straightforward exercise to show that, as $\lambda\to+\infty$
along the real axis, \eqnref{frozenchar} has two kinematic solutions
$$
\mu\sim  \pm \lambda^{1/2},
\quad \pmatrix
U \\ Z
\endpmatrix
\sim \pmatrix 
1\\ 0
\endpmatrix,
$$
or
$
(u, u', z)
\sim e^{\mu x} 
(1, \pm \lambda^{1/2}, 0)$
in phase plane variables,
and one reactive solution
$$
\mu\sim  \lambda/s,
\quad \pmatrix
U \\ Z
\endpmatrix
\sim \pmatrix 
0\\ 1
\endpmatrix,
$$
or 
$
(u, u', z)
\sim e^{\mu x} 
(0, 0, 1)$ in phase variables.
Thus, the stable and unstable subspaces of the frozen-coefficient
equations have a uniform spectral gap of order $\lambda^{1/2}$,
and, carefully applying the tracking lemma of [GZ,ZH,Z.3], we may 
conclude that the stable and unstable manifolds of the variable-coefficient
equations \eqnref{15} track the corresponding subspaces of the 
frozen-coefficient equations to angle $\lambda^{-1/2}$.
See [Z.3] and especially [MZ.3] for similar calculations.

>From these considerations, we find that
$$
\pmatrix
u_1^+\\u_1^{+'} \\ z_1^+
\endpmatrix
\sim\pmatrix 
1\\ \sqrt{\lambda}\\ 0
\endpmatrix
\eta_+(x),
$$
$$
\pmatrix
u_2^- & u_3^- \\ 
u_2^{-'} & u_3^{-'} \\ 
z_2^- & z_3^-
\endpmatrix
\sim
\pmatrix
0 & 1\\
0 & -\sqrt{\lambda} \\
1 & 0
\endpmatrix
\eta_-(x),
$$
where $\eta_+(x)$ is a real nonvanishing scalar function,
and $\eta_-(x)$ is a real nonsingular $2\times 2$ matrix function.
Thus,
$$
\eqalign{
D &\sim \det
\pmatrix
1&0&1\\
\sqrt{\lambda }&0&-\sqrt{\lambda }\\
0&\lambda &0
\endpmatrix 
\eta_+ \det \eta_- \cr
&= \eta_+ \det \eta_- \cr
&= u_1^+ \det\
\pmatrix
u_2^- & u_3^- \\
z_2^- & z_3^-
\endpmatrix,
}
$$
from which we may deduce the claim by reality of $D$ and
nonvanishing of $\eta_+$, $\det \eta_-$.

Alternatively, following [BSZ], we may first deduce by
by standard G\"arding-type energy estimates
that there exist no zeroes of $D$/eigenvalues of $L$ 
for sufficiently large real $\lambda$, uniformly, 
independent of (bounded) model parameters.
Performing the homotopy $f$, $k\to0$ reduces the eigenvalue
equations to the simple case
$$
\cases
u''=\lambda u,\\
z'=(\lambda/s)z,
\endcases
$$
for which the claim follows by an elementary direct calculation.
\myqed
\medskip

{\bf The stability index.} Combining our calculations at $0$ and $+\infty$,
we have

\proclaim{Corollary \thmlbl{strongindex}}
For nondegenerate strong detonations, the stability index 
$\Gamma:=\sgn D'(0)D(+\infty)$ is well-defined,
and satisfies
$$
\Gamma= \sgn([u]+q) \gamma^2 \bar u_x(+\infty)>0.
\eqnlbl{31}
$$
Thus, the number of zeroes of $D$ on the unstable half-plane is
even, consistent with stability.
\endproclaim

{\bf Proof.} $D(\lambda)$ is real for real $\lambda$, and, by 
the previous lemma, nonvanishing for real $\lambda$ sufficiently
large.  Thus, $\sgn D(\lambda)$ is independent of $\lambda$ for
real $\lambda$ sufficiently large, and $\sgn D(+\infty)$ is well-defined.
Moreover, since $D$ is analytic, the number of positive
real roots of $D$ is even or odd according as $\Gamma$ is positive or negative
(indeterminate parity if $\Gamma=0$, corresponding to $D'(0)=0$). 
In fact, since nonreal roots occur in conjugate pairs, 
by complex symmetry of $D$, the sign of $\Gamma$ determines the parity
of the number of zeroes of $D$ on the entire unstable half-plane
$\{\lambda:\, \R \lambda>0\}$.
Thus, it is sufficient to verify \eqnref{31} in order to establish
the claim.

To establish the formula $\Gamma= \sgn([u]+q) \gamma^2 \bar u_x(+\infty)$,
it is sufficient to show that expression \eqnref{largelambda} is independent
of $\lambda$ for real $\lambda\ge 0$;
evaluating at $\lambda=0$ and combining with the formula
of Proposition \thmref{4} for $D'(0)$, we then obtain the result.
Here, we are using also the fact that 
$$\sgn \gamma=\sgn \det \pmatrix u^-_2 & u_3^-\\ 
 z^-_2& z_3^-
\endpmatrix (0)
=
\sgn \det \pmatrix u^-_2 & u_3^-\\ 
 z^-_2& z_3^-
\endpmatrix (-\infty)
\ne 0
$$
as a nontrivial (since its columns are
by construction independent as $x\to -\infty$) Wronskian for the linearized
traveling-wave ODE.

To establish independence of \eqnref{largelambda} with respect to $\lambda$,
we may prove a projection lemma as in [GZ], namely that projection
onto $u$ coordinate is full rank on the unstable subspace, and projection
onto $(u,z)$ coordinates is full rank on the stable subspace of the coefficient
matrix of the limiting eigenvalue equation.  This is straightforward in
the present case,
using the limiting structure $(1, \mu, 0)^t$ of unstable eigenvectors, and
$$
\pmatrix
1&*\\
\mu & * \\
0 &1
\endpmatrix
$$
of stable eigenvectors at $x\to \pm \infty$.

Finally, the conclusion 
$\sgn([u]+q) \gamma^2 \bar u_x(+\infty)>0$
follows from the facts that:
(i) $\bar u_x(+\infty)<0$ for all strong
detonation profiles, by our earlier phase plane analysis, 
(ii) $[u]+q=(1/s)[f(u)]<0$ by $[u]<0$, $s>0$, 
and the fact that $f$ is strictly monotone increasing,
and (iii) $\gamma \ne 0$. 
\myqed

{\bf Remark \thmlbl{ext}.} A similar result has been shown in the ZND limit
($k$ sufficiently small) in [Ly,LyZ],  
for the full, Navier--Stokes equations of reacting flow with
an ideal gas equation of state.
Multi-dimensional analogues are given in [Z.3,JLy.1--2].
\smallskip

{\bf Remark \thmlbl{hist}.} 
A determinant analogous to the 
Evans function was used by Erpenbeck [Er.2,Er.4] 
to study stability of discontinuous strong detonation
fronts within the framework of the ZND model; 
moreover, he obtained an interesting stability criterion 
by examining the {\it high} frequency limit, opposite to what we do
(and outside the low-frequency regime in which the ZND model
is expected to describe model \eqnref{1abstract}) [Er.3];
for further discussion, see [Z.3], Appendix A.3, or [JLy.2].  
The stability criteria of Erpenbeck--Majda for
shock waves [Er.1,M.1--3], of Erpenbeck for detonations [Er.1], 
and the Evans function criterion of Evans--Alexander--Gardner--Jones for
reaction--diffusion equations [E.1--4,J.1,AGJ] are thus seen to be 
facets of the same general principle/stability criterion, applied
to equations of varying type.  For further discussion, see [ZS].

\bigskip
\newsection {\bf Section 7. Stability analysis for weak detonations.}
\sectionnumber=7
\theoremnumber=0
\equationnumber=0
\smallskip
\TagsOnLeft

We next consider the interesting weak detonation case
$\alpha _-$, $\alpha _+<0$, for which the
speed of the detonation everywhere exceeds the rate of propagation of gas 
dynamical signals.  These correspond to {\it undercompressive}, saddle--saddle
connections and connect $u_+$ only to special states $u_-$, as
discussed in Section 3.
As usual, we assume nondegeneracy, i.e., spatially exponential decay.

Repeating the stability analysis of the previous section, 
we now find that, at $\lambda =0$,
the mode $(u^-_2,z^-_2)$ is not asymptotically decaying, but
asymptotically {\it constant},
$$
\pmatrix
u^-_2\\
{u^-_2}'\\
z^-_2
\endpmatrix\quad \longrightarrow \quad
\pmatrix
1\\
0\\
0
\endpmatrix \eqnlbl{32}
$$
as $x\to -\infty$
(see \eqnref{17}).  Thus for $j=2$, we obtain in
place of \eqnref{25} the equation (recall, $\alpha_-<0$)
$$
{u^-_2}' - \alpha u^-_2 + sqz^-_2=-\alpha _-=|\alpha_-|.
\eqnlbl{33}
$$

\proclaim{Proposition \thmlbl{6}}  For nondegenerate weak detonations, 
$D(0)=0$,
while
$$
D'(0)=(\partial d/\partial s) (u_+,s) |\alpha _-|<0,
\eqnlbl{33.1}
$$
where $d(u_+,s)$ is the Melnikov separation function 
defined in \eqnref{melnikov},
$\partial d/\partial s$ is as given in  \eqnref{33.2},
and $\alpha$ as usual denotes $a-s= df(\bar u) -s$.
\endproclaim

{\bf Proof}.  Following the steps of the previous case, we obtain in
place of \eqnref{30}:
$$
D'(0)=\det
\pmatrix
\barU_x & u^-_2 &y\\
0&|\alpha _-| & [u]\\
\Zbar_x &z^-_3 & \sigma
\endpmatrix_{|_{x=0}},
\eqnlbl{34}
$$
$y,\sigma$ as defined previously.  
(Note: the only difference is in the new 2,2 entry $|\alpha_-|$
coming from the righthand side of \eqnref{33}.)
This in turn gives
$$
\aligned
D'(0)
&=\det
\pmatrix
\barU_x & u^-_2 &\tildeY\\
0 &|\alpha _-| & 0\\
\Zbar_x &z^-_3 & \Tsigma
\endpmatrix_{|_{x=0}}\\
&=|\alpha _-| \det
\pmatrix
\barU_x &\tildeY\\
\Zbar_x & \Tsigma
\endpmatrix_{|_{x=0}},
\endaligned
\eqnlbl{35}
$$
where $\tildeY:=y-\Big(\frac{[u]}{|\alpha _-|}\Big)u^-_2$,
$\Tsigma:=\sigma - \Big(\frac{[u]}{|\alpha _-|}\Big)z^-_2$.  But, 
$I=\det 
\pmatrix 
\barU_x & \tildeY\\
\Zbar_x & \Tsigma
\endpmatrix_{|_{x=0}}$, just as in the standard, undercompressive shock
case treated in [GZ], 
can be recognized as the Melnikov integral \eqnref{33.2},
by expressing $I=I^--I^+$, where
$$
I^\pm :=
\pmatrix 
\barU_x & \tildeY^\pm\\
\Zbar_x & \Tsigma^\pm
\endpmatrix_{|_{x=0}},
$$
with $(\tilde y^+,\tilde \sigma^+)=(y^+,\sigma^+)$
and 
$(\tilde y^-,\tilde \sigma^-)=(y^-,\sigma^-)  
-([u]/|\alpha _-|)(u^-_2,z_2^-)$,
and observing that $I^\pm$ (since do $\tilde y^\pm$, $\tilde \sigma^\pm$)
satisfy the same ODE
$$
I'-(\alpha + k\varphi)I = \det
\pmatrix 
\barU_x & (\bar u- u_+) + q(\bar z-z_+)\\
\Zbar_x & (\bar z_x/s)
\endpmatrix_{|_{x=0}},
$$
and using Duhamel's principle to express $I^\pm$ as the
integral from $\pm \infty$ to $0$ of the integrand in the first
line of \eqnref{33.2};
see, e.g., [GZ] for details. 
The evaluation in the second line follows from the identity
$-s(\bar u-u_+)-sq(\bar z-z_+)= \bar u'- (f(\bar u)-f(u_+))$,
a rearrangement of the traveling-wave ODE, and elementary calculation.
The sign of the integral then follows from $\bar z_x>0$, 
$\bar u_x<0$, and $f'>0$.
\myqed

\proclaim{Corollary \thmlbl{fullweak}}
For nondegenerate weak profiles, there holds
$$
\Gamma=\sgn |\alpha_-| (\partial d/\partial(u_+,s))
\bar u_x(+\infty) \bar z_x(-\infty) >0,
\eqnlbl{wdstable}
$$
consistent with stability.
\endproclaim

{\bf Proof.}  The result of Proposition \thmref{5} is
independent of the case, and so the evaluation \eqnref{largelambda} of 
$\sgn D(+\infty)$ remains valid, as does the homotopy to $\lambda=0$.
However, $\det \pmatrix u_2^- & u_3^-\\ z_2^- & z_3^-\endpmatrix$
now becomes asymptotic to the explicitly prescribed endstates
$$
\det \pmatrix 1 & \bar u_x\\ 0 & \bar z_x \endpmatrix(-\infty)
= \bar z_x(-\infty),
$$
giving a combined result of 
$\sgn D(+\infty)=\sgn \bar u_x(+\infty)\bar z_x(-\infty)<0$.  
Combining with our previous calculation of $\sgn D'(0)<0$,
we obtain the result.
\myqed

{\bf Remark \thmlbl{wellp}.} As in [GZ], notice the relation between the key 
Melnikov integral $\partial d/\partial s$ 
and well-posedness of the Riemann problem,
Remark \thmref{linwp}.

\bigskip
\newsection {\bf Section 8. Stability of degenerate waves.}
\sectionnumber=8
\theoremnumber=0
\equationnumber=0
\smallskip
\TagsOnLeft

Our methods, in the basic form presented here, 
require spatially exponential decay of the profile,
so do not apply to the remaining class of degenerate,
subalgebraically decaying profiles. \footnote{
See however [H.1--2,HZ,SS] for extensions to more general situations.}
Nonetheless, we can make one or two comments of a general nature
regarding their stability.
\medskip
{\it Degenerate detonation profiles.}
Degenerate detonation profiles occur only in the special case
$u_+=u_i$, and generically do not appear in Riemann solutions;
in this sense, they may be considered a technical curiosity 
arising from the physically artificial cutoff at $u_i$ of the ignition
function (for further discussion, see Section 9).
When they do appear, they occur as a one-parameter family of
degenerate strong detonation profiles, bounded above by
a nondegenerate, and presumably stable, strong detonation profile,
and below by a nondegenerate gas-dynamical shock profile followed by
a degenerate weak detonation profile.
%
%

This is somewhat reminiscent of the one-parameter families of overcompressive
shocks studied for nonstrictly hyperbolic conservation laws
in [Fre.2,L,FreL].
As in that situation, degenerate strong detonations cannot be 
singly orbitally stable, but only as a family; 
moreover, in the small viscosity limit (i.e., the true ZND limit, without
rescaling), there is the same situation that the $L^1$ difference between
profiles corresponding to fixed orbits goes to zero, so that stability
if it exists at all cannot be expected to be uniform with respect to
$L^1$ in the small viscosity limit.
(Note that, since subalgebraically
decaying, they are infinitely far in $L^1$ from the bounding, 
exponentially decaying strong detonation.)
On the other hand, conservation of mass (in this case, of $\int (u+qz)dx$)
no longer determines uniquely the time-asymptotic profile, and
it is not clear that any degenerate profile should be expected to
be stable under perturbation even for fixed viscosity.

In any case, at the level of Riemann solutions, 
degenerate strong detonations cannot be distinguished from the 
bounding strong detonation,
so give no contradiction with the conclusions of [CF] based on the ZND
model. 
On the other hand, as we saw in Section 4, existence of degenerate 
weak detonation profiles does lead to nonuniqueness in solutions of the Riemann
problem for the special case $u_+=u_i$.  
One might hope, therefore, that they could be discarded on 
grounds of instability, thereby validating ZND conclusions for 
$k$ sufficiently small.

\medskip
{\it Degenerate deflagration profiles.}
Contrary to the situation of degenerate detonations,
weak deflagrations are always of degenerate type, and
play an important role in the solution of Riemann problems.
In particular, stability of (degenerate) weak deflagration profiles is
necessary for any $k>0$ in order that there exist a general 
Riemann solution composed of stable viscous profiles, hence 
expected from the ZND point of view.
Thus, stability vs. instability of these waves has important
philosophical consequences, and is deserving of further study.
We point out only that, discarding cross terms in the linearized
equations, all remaining terms for weak deflagrations are favorable for
a basic $L^2$ energy estimate, in contrast to the 
detonation case: the term $(\alpha u)_x$,
because $\alpha_x \ge 0$ (expansivity) due to monotonicity 
of the wave, 
and the term $qk\bar d\varphi \bar z u$ 
because $d\varphi \le 0$ in the deflagration regime.

The stability of degenerate waves of either type (detonation
or deflagration) is likely to be sensitive due to the subexponential
rate of decay of the background profile,
as for example in the analogous case of KPP waves; see, e.g., [He].

\bigskip
\newsection {\bf Section 9. Extensions.}
\sectionnumber=9
\theoremnumber=0
\equationnumber=0
\smallskip
\TagsOnLeft

Finally, we discuss various elaborations that can be
accomodated in the basic model without affecting the analysis:

\medskip
{\bf Multi-species and multiple reactions}.
Combustion involving more than one reactant or reaction can be modeled
abstractly, following [FD], by the use of {\it progress variables}
$$
\Lambda =(\lambda _1,\lambda _2,\cdots,\lambda _m),
\eqnlbl{201}
$$
each $0\le \lambda _j \le 1$ denoting the progress toward completion
of a single idealized reaction.  As described in [FD], the same
framework can be used to describe arbitrarily complex reactions in a
compact form,with the $\lambda _j$ now representing linear
combinations of progress variables for several
simpler reactions.  The progress variables satisfy a {\it rate
equation} 
$$
\dot{\Lambda }=R(\Lambda ,u),
\eqnlbl{202}
$$
modeling the composite chemical reactions, and enter the kinematic
equations via
$$
(u- Q\Lambda)_t+(f(u)=u_{xx},
\eqnlbl{203}
$$
where $Q=(q_1,q_2,\cdots q_m)$ denote heats of reaction, $q_j>0$ for
{\it exothermic reactions}, for $q_j<0$ for {\it endothermic
reactions}. 

For comparison, consider first a single reaction $A\to B$.  This is
represented in the present notation by a single progress variable
$\lambda _1$, with rate equation
$$
\dot{\lambda}_1 = k(\lambda _1,u)(1-\lambda _1),
\eqnlbl{204}
$$
a typical physical rate $k$ being the Arrhenius rate 
$$
k(\lambda _1,u)=e^{-{E_0 \over RT}},
\eqnlbl{205}
$$
where $T$ is temperature, $R$ is the gas constant, and $E_0$ is 
{\it activation energy}: the threshold determining ignition temperature. 
In the simple model considered previously, $z$ corresponds
to $(1-\lambda _1)$, and $k\varphi$ to the rate function $k (\cdot,
\cdot)$. 
Modification of the Arrhenius rate by a low-temperature
cutoff $T_i$ is a standard device in the combustion literature, circumventing
the ``cold boundary'' problem that the Arrhenius kinetics do
not permit a stationary unburned state and therefore preclude the
existence of traveling wave profiles.

A more realistic model of combustion consisting of two, consecutive
reactions, $A\to B, B\to C$, is modeled by the rate equations
$$
\cases
\dot{\lambda }_1 = k_1(\lambda,u)(1-\lambda _1),\\
\dot{\lambda }_2 = k_2(\lambda,u)(\lambda _1-\lambda _2),
\endcases
\eqnlbl{206}
$$
with either (i)  $q_1$, $q_2>0$ (exothermic/exothermic), or 
(ii)  $q_1 > 0$, $q_2 < 0$ (exothermic/endothermic).
As above, a standard choice of rate function would be
$k_j=\varphi_j(u)k_j$, with $k_j$ constant, and $\varphi_j$
ignition functions of the usual form.
%
As pointed out in [FD], the traveling wave structure in case
(i) is very similar to that for a single reaction; on the other hand,
(ii) introduces important new phenomena at the ZND level, 
explaining e.g. an observed experimental shift from CJ to weak 
detonation in ignition problem.
See [FD], p. 168--173, for further discussion.  
Three or more reactions can easily be
included, but do not appear to result in significant new phenomena
[FD]. 

It is straightforward to show that the main conclusions of the previous 
sections remain valid for equations \eqnref{206}, or indeed any system of
reactions having such a triangular structure: 
more generally, for any reaction system in which 
$\Lambda =(1,\cdots1)$ is stable above the ignition temperature(s) 
as a rest point of rate equation \eqnref{206}. 
In particular, we recover the result for strong detonations that
the stability index is always positive, consistent with stability.
Likewise, we recover the same expression for $\Gamma$ in the case
of weak detonations, but now the Melnikov integral for
$(\partial d/\partial s)_{(u_+,s)}$ is more complicated and its
sign (apparently) no longer explicitly evaluable.
It would be very interesting to investigate structure and stability of 
traveling waves for the 2--reaction model above, especially in light
of anomalous behavior (e.g., a shift in the ignition problem
from CJ to weak detonation) predicted in [FD,B] in the the ZND setting,
in the high--activation energy limit.

\medskip
{\bf Reaction-dependent equation of state}.
As pointed out in [CHT], a more realistic assumption is
that the gas-dynamical equation of state (EOS) $f$
depend not only on $u$, but also on $z$,
through the chemical makeup of the gas.
With this change, the linearized equations \eqnref{13}, \eqnref{13alt}
become
$$
\cases
u_t 
-q( k\varphi '(\barU)u \Zbar - k \varphi(\barU)z)
+ (\alpha  u)_x + (\beta z)_x  = u_{xx}\\
z_t - sz_x = 
-\ k\varphi '(\barU)u \Zbar -
 k \varphi(\barU)z,
\endcases\eqnlbl{13rd}
$$
$$
\cases
u_t + qz_t + (\alpha  u)_x+(\beta z)_x - sqz_x = u_{xx}\\
z_t - sz_x = 
-\ k\varphi '(\barU)u \Zbar - k \varphi(\barU)z,
\endcases\eqnlbl{13}
$$
where $\alpha:= f_u(\bar u, \bar z)(x)-s$ and $\beta:= f_v(\bar u, \bar z)(x)$.
It is straightforward to check that all of our stability analysis goes
through as before in this more general setting, 
to yield exactly the same geometric stability conditions.
In particular, $\Gamma>0$
holds always for a strong detonation profile if it exists.
Likewise, the expression derived for $\Gamma$ in the weak
detonation case, including our calculation of the Melnikov integral 
$\partial d/\partial s$, remains valid,  though the sign
of the integral may change depending on the details of the equation of state.

Regarding the existence problem, and the related question of stability
of weak detonation profiles, we are led naturally to a simple condition,
$$
f_z \le 0,
\eqnlbl{monrd}
$$
under which all of the conclusions of the paper generalize to the
case of a reaction-dependent EOS.
For, this implies that $F_z=f_z-qs<0$, where $F=0$ as in \eqnref{F}
defines the nullcline $u'=0$ for the traveling-wave equations,
and so the nullcline is again a graph over $z$.
This was essentially the only property used in the qualitative analysis
of the phase portrait, and so we recover all related conclusions,
including the important property of monotonicity in $u$ and $z$ 
of weak detonation profiles.
Likewise, \eqnref{monrd} implies that 
$f(\bar u,\bar z)-f(u_+,z_+)>0$ along weak detonation profiles,
since $f_u>0$, $f_z\le 0$ and $\bar u$ and $\bar z$ are, respectively,
monotone decreasing and increasing.
Thus, we obtain the key monotonicity in $s$ asserted in
Proposition \thmref{mon}, now in the restricted case $\hat z\le 1$,
in particular, the conclusion that $\partial d/\partial s<0$ at
a profile $d=0$, along with the other stated monotonicity results
(the arguments for which are unaffected by dependence of $f$ on $z$),
and thereby the remaining conclusions of the paper.

On the other hand, when the monotonicity condition \eqnref{monrd} 
is violated, we see no reason why anomalies in both existence and 
stability theory might not occur, and this seems an important
direction for further investigation.
In the full, reacting Navier--Stokes equations, the flux functions,
taking gas composition into account, are just the convex averages,
weighted by $z$, of corresponding fluxes for pure unburned and pure burned
gas, and similarly for the coefficient $c_v$ relating temperature
to internal energy.
We conjecture that the corresponding monotonicity condition in 
this physical setting reduces to the condition that the gas constant
$\gamma(1)$ of the unburned gas be less than the gas constant
$\gamma(0)$ of the burned gas, assuming that each separately obeys 
the EOS of an ideal,  polytropic gas: i.e., from the 
kinetic theory of gases point of view,
the average number of internal degrees of freedom per molecule $n$ 
{\it decrease}
upon reaction (recall: $\gamma=(n+2)/n$; see, e.g., [Ba], pp. 37--45).
The derivation of an analog of condition \eqnref{monrd} by analysis
of the ZND limit we regard as another very interesting open problem.

\bigskip
\Refs
\medskip\noindent
[AT] G. Abouseif and T.Y. Toong,
{\it Theory of unstable one-dimensional detonations,}
Combust. Flame 45 (1982) 67--94.
\medskip\noindent
[AGJ] J. Alexander, R. Gardner and C.K.R.T. Jones,
{\it A topological invariant arising in the analysis of
traveling waves}, J. Reine Angew. Math. 410 (1990) 167--212.
\medskip\noindent
[AlT] R. L. Alpert and T.Y. Toong,
{\it Periodicity in exothermic hypersonic flows
about blunt projectiles,} Acta Astron. 17 (1972) 538--560.
\medskip\noindent\noindent
[AMPZ.1] A. Azevedo, D. Marchesin, B. Plohr, and K. Zumbrun,
{\it Nonuniqueness of solutions of Riemann problems,}
 Z. Angew. Math. Phys. 47 (1996) 977--998. 
\medskip\noindent
[BE] A.A. Barmin and S.A. Egorushkin,
{\it Stability of shock waves,}
Adv. Mech. 15 (1992) No. 1--2, 3--37.
\medskip\noindent
[Ba] G.K. Batchelor,
{\it An introduction to fluid dynamics,} 
Second paperback edition. Cambridge Mathematical
Library. Cambridge University Press, Cambridge (1999) xviii+615 
pp. ISBN: 0-521-66396-2.
\medskip\noindent
[BSZ] S. Benzoni--Gavage, D. Serre, and K. Zumbrun,
{\it Alternate Evans functions and viscous shock waves,}  
to appear, SIAM J. Math. Anal. 
\medskip\noindent
\medskip\noindent
[BMR] A. Bourlioux, A. Majda, and V. Roytburd,
{\it Theoretical and numerical structure for unstable one-dimensional
detonations,} SIAM J. Appl. Math. 51 (1991) 303--343.
\medskip\noindent
[Br.1] L. Q. Brin, {\it Numerical testing of the stability of viscous
shock waves}, Ph.D. dissertation, Indiana University, May 1998.
\medskip\noindent
[Br.2] L. Q. Brin,
{\it Numerical testing of the stability of viscous shock waves,} 
Math. Comp. 70 (2001), no. 235, 1071--1088.
\medskip\noindent
[BrZ] L. Brin and K. Zumbrun, 
{\it Analytically varying eigenvectors and the stability of viscous
shock waves,} to appear, Mat. Contemp. (2003).
\medskip\noindent
[B] J.D. Buckmaster, {\it An introduction to combustion theory,}
3--46,
in {\it The mathematics of combustion,} Frontiers in App. Math.
(1985) SIAM, Philadelphia ISBN: 0-89871-053-7.
\medskip\noindent
[BL] J.D. Buckmaster and G.S.S. Ludford, 
{\it The effect of structure on the
stability of detonations I. Role of the induction zone,}
1669--1675, 
in {\it Proc. Twenty-First Symp. (Intl) on Combustion} (1988)
The Combustion Institute, Pittsurgh, PA.
\medskip\noindent
[BN] J.D. Buckmaster and J. Nevis, 
{\it  One-dimensional detonation stability- The spectrum for
infinite activation energy,}
Phys. Fluids 31 (1988) 3571--3575.
\medskip\noindent
[B] B. Bukiet, 
{\it The effect of curvature on detonation speed,}
SIAM J. Appl. Math. 49 (1989) 1433--1446.
\medskip\noindent
[CHT] G.-Q. Chen, D. Hoff, and K. Trivisa,
{\it Global solutions to a model for
exothermically reacting, compressible flows with large discontinuous
initial data,}  preprint (2002).
\medskip\noindent
[CF] R. Courant and K.O. Friedrichs,
{\it Supersonic flow and shock waves,}
Springer--Verlag, New York (1976) xvi+464 pp. 
\medskip\noindent
[CMR] P. Collella, A. Majda, and V. Roytburd,
{\it Theoretical and numerical structure for reacting shock waves,}
SIAM J. Sci. Stat. Comput. 7u (1986) 1059--1080.
\medskip\noindent
[Er.1] J. J. Erpenbeck,
{\it Stability of step shocks,} Phys. Fluids 5 (1962) no. 10, 1181--1187.
\medskip\noindent
[Er.2] J. J. Erpenbeck,
{\it Stability of steady-state equilibrium detonations,}
Phys. Fluids 5 (1962) no. 5, 1181--1187.
\medskip\noindent
[Er.3] J. J. Erpenbeck, 
{\it Detonation stability for disturbances
of small transverse wavelength,}
Phys. Fluids 9 (1966) No. 7,1293--1306.
\medskip\noindent
[Er.4] J. J. Erpenbeck, 
{\it Stability of idealized one--reaction detonations,}
Phys. Fluids 7 (1964) No. 5, 684--696.
\medskip\noindent
[Er.5] J. J. Erpenbeck, 
{\it Steady detonations in idealized two--reaction systems,}
Phys. Fluids 7 (1964) No. 9, 1424--1432.
\medskip\noindent
[Er.6] J. J. Erpenbeck, 
{\it Nonlinear theory of unstable
one--dimensional detonations,}
Phys. Fluids 10 (1967) No. 2,274--289.
\medskip\noindent
[Er.7] J. J. Erpenbeck, 
{\it Structure and stability of the square-wave detonation,}
in {\it Ninth Symp. (intl) on Combustion} (1963) 442--453.
\medskip\noindent
[E.1] J.W. Evans,
{\it Nerve axon equations: I. Linear approximations,}
Ind. Univ. Math. J. 21 (1972) 877--885.
\medskip\noindent
[E.2] J.W. Evans,
{\it Nerve axon equations: II. Stability at rest,}
Ind. Univ. Math. J. 22 (1972) 75--90.
\medskip\noindent
[E.3] J.W. Evans,
{\it Nerve axon equations: III. Stability of the nerve impulse,}
Ind. Univ. Math. J. 22 (1972) 577--593.
\medskip\noindent
[E.4] J.W. Evans,
{\it Nerve axon equations: IV. The stable and the unstable impulse,}
Ind. Univ. Math. J. 24 (1975) 1169--1190.
\medskip\noindent
[F.1] W. Fickett,
{\it Stability of the square wave detonation in a model
system,}  Physica 16D (1985) 358--370.
\medskip\noindent
[F.2] W. Fickett, {\it Detonation in miniature,} 133--182,
in {\it The mathematics of combustion,} Frontiers in App. Math.
(1985) SIAM, Philadelphia ISBN: 0-89871-053-7.
\medskip\noindent
[FD] W. Fickett and W.C. Davis,
{\it Detonation,} University of California Press, Berkeley, CA (1979):
reissued as {\it Detonation: Theory and experiment,}
Dover Press, (2000), ISBN 0-486-41456-6.
\medskip\noindent
[FW] Fickett and Wood, 
{\it Flow calculations for pulsating one-dimensional
detonations,} Phys. Fluids 9 (1966) 903--916. 
\medskip\noindent
[Fo] G.R. Fowles,
{\it On the evolutionary condition for stationary plane
waves in inert and reactive substances,}
in {\it  Shock induced transitions and phase structures in general media,}
\medskip\noindent
[Fre.1] H. Freist\"uhler,
{\it A short note on the persistence of ideal shock waves,}
Arch. Math. (Basel) 64 (1995) 344--352.
\medskip\noindent
[Fre.2] H. Freist\"uhler,
{\it Dynamical stability and vanishing viscosity:
A case study of a nonstrictly hyperbolic system of conservation laws},
Comm. Pure Appl. Math. 45 (1992) 561--582.
\medskip\noindent
[FreL] H. Freist\"uhler and T.-P. Liu,
{\it Nonlinear stability of overcompressive shock
waves in a rotationally invariant system of viscous conservation laws},
Commun. Math. Phys. 153 (1993) 147--158.
\medskip\noindent
[FreS] H. Freist\"uhler and P. Szmolyan, 
{\it Spectral stability of small shock waves,}
Arch. Ration. Mech. Anal. 164 (2002) 287--309. 
\medskip\noindent
[GJ.1] R. Gardner and C.K.R.T. Jones,
{\it A stability index for steady state solutions of
boundary value problems for parabolic systems},
J. Diff. Eqs. 91 (1991), no. 2, 181--203. 
\medskip\noindent
[GJ.2] R. Gardner and C.K.R.T. Jones,
{\it Traveling waves of a perturbed diffusion equation
arising in a phase field model}, 
Ind. Univ. Math. J. 38 (1989), no. 4, 1197--1222.
\medskip\noindent
[G] R.A. Gardner, 
{\it On the detonation of a combustible gas,}
Trans. Amer. Math. Soc. 277 (1983) 431--468.
\medskip\noindent
[GS.1] I. Gasser and P. Szmolyan,
{\it A geometric singular perturbation analysis of detonation 
and deflagration waves,} SIAM J. Math. Anal. 24 (1993) 968--986.
\medskip\noindent
[GS.2] I. Gasser and P. Szmolyan,
{\it Detonation and deflagration waves with multistep reaction schemes,}
SIAM J. Appl.  Math. 55 (1995) 175--191. 
\medskip\noindent
[GZ] R. Gardner and K. Zumbrun,
{\it The Gap Lemma and geometric criteria for instability
of viscous shock profiles}, 
Comm. Pure Appl.  Math. 51 (1998), no. 7, 797--855. 
\medskip\noindent
[JLy.1] K. Jenssen and G. Lyng, {\it Multidimensional stability
of viscous detonation waves,} in preparation.
\medskip\noindent
[JLy.2] K. Jenssen and G. Lyng, {\it Multidimensional stability
of ZND detonation waves,} in preparation.
\medskip\noindent
[GH] J. Guckenheimer and P. Holmes,
{\it Nonlinear oscillations, dynamical systems, and bifurcations of vector fields},
(Revised and corrected reprint of the 1983 original),
Springer--Verlag, New York (1990), xvi+459 pp.
\medskip\noindent
[HK] J. Hale and H. Kocak,
{\it Dynamics and bifurcations,}
Texts in Applied Mathematics, 3. Springer-Verlag, New York, 1991.
xiv+568 pp. ISBN: 0-387-97141-6. 
\medskip\noindent
[He] D. Henry,
{\it Geometric theory of semilinear parabolic equations},
Lecture Notes in Mathematics, Springer--Verlag, Berlin (1981),
iv + 348 pp.
\medskip\noindent
[H.1] P. Howard,
{\it Pointwise estimates and stability for degenerate viscous shock waves,}
J. Reine Angew. Math.  545 (2002) 19--65. 
\medskip\noindent
[H.2] P. Howard,
{\it Local tracking and stability for degenerate viscous shock waves,}
J. Differential Equations 186 (2002) 440--469. 
\medskip\noindent
[HZ] P. Howard and K. Zumbrun,
{\it The Evans function and 
stability criteria for degenerate viscous shock waves,}
preprint (2002).
\medskip\noindent
[J] C.K.R.T. Jones,
{\it Stability of the travelling wave solution of the FitzHugh--Nagumo system},
Trans. Amer. Math. Soc.  286 (1984), no. 2, 431--469.
\medskip\noindent
[KS] T. Kapitula and B. Sandstede,
{\it Stability of bright solitary-wave solutions 
to perturbed nonlinear Schrödinger equations,} Phys. D 124
(1998), no. 1-3, 58--103.
\medskip\noindent
[Li.1] T. Li,
{\it On the Riemann problem for a combustion model,}
SIAM J. Math.  Anal. 24 (1993), no. 1, 59--75.
\medskip\noindent
[Li.2] T. Li,
{\it On the initiation problem for a combustion model,}
J. Differential Equations 112 (1994), no. 2, 351--373.
\medskip\noindent
[Li.3] T. Li,
{\it Rigorous asymptotic stability of a Chapman--Jouguet detonation
wave in the limit of small resolved heat release,}
Combust. Theory Model. 1 (1997), no. 3, 259--270.
\medskip\noindent
[Li.4] T. Li,
{\it Time-asymptotic limit of solutions of a combustion problem,}
J.  Dynam. Differential Equations 10 (1998), no. 4, 577--604.
\medskip\noindent
[Li.5] T. Li,
{\it Stability of strong detonation waves and rates of convergence,}
Electron. J. Differential Equations (1998) no. 9, 17 pp. (electronic).
\medskip\noindent
[Li.6] T. Li, Zurich proc., 
{\it Stability and instability of detonation waves,}
Hyperbolic problems: theory, numerics, applications, 
Vol. II (Z\"{u}rich, 1998), 641--650, Internat. Ser. Numer. 
Math., 130, Birkhäuser, Basel, 1999.
\medskip\noindent
[Li.7] T. Li, 
{\it Stability of a transonic profile arising from divergent 
detonations,}  Comm. Partial Differential Equations 25 (2000) 
11-12, 2087--2105.
\medskip\noindent
[LS] H.I. Lee and D.S. Stewart,
{\it Calculation of linear detonation instability:
one-dimensional instability of plane detonation,}
J. Fluid Mech. (1990) 102--132.
\medskip\noindent
[L] T.-P. Liu,
{\it Nonlinear stability and instability of overcompressive shock
waves,} in:
{\it Shock induced transitions and phase structures in general media},
159--167, IMA Vol.  Math. Appl., 52, Springer, New York, 1993.
\medskip\noindent
[LLT] D. Li, T.-P. Liu, and D. Tan,
{\it Stability of strong detonation travelling waves to combustion model,}
J. Math. Anal. Appl. 201 (1996), no. 2, 516--531.
\medskip\noindent
[LYi] T.-P. Liu and L.A. Ying,
{\it Nonlinear stability of strong detonations for
a viscous combustion model,}
SIAM J. Math. Anal. 26 (1995), no. 3, 519--528.
\medskip\noindent
[LY] T.-P. Liu and S.-H. Yu,
{\it Nonlinear stability of weak detonation waves for a combustion model,}
Comm. Math. Phys. 204 (1999), no. 3, 551--586. 
\medskip\noindent
[LZ.1] T.P. Liu and K. Zumbrun,
{\it Nonlinear stability of an undercompressive shock for complex
Burgers equation,} Comm. Math. Phys. 168 (1995), no. 1, 163--186.
\medskip\noindent
[LZ.2] T.P. Liu and K. Zumbrun,
{\it On nonlinear stability of general undercompressive viscous shock waves,}
Comm.  Math. Phys. 174 (1995), no. 2, 319--345.
\medskip\noindent
[Ly] G. Lyng, {\it One dimensional stability of detonation waves,}
doctoral thesis, Indiana University (2002).
\medskip\noindent
[LyZ] G. Lyng and K. Zumbrun,  {\it Stability of 
detonation waves,} preprint (2003).
\medskip\noindent
[M.1], A. Majda,
{\it The stability of multi-dimensional shock fronts -- a
new problem for linear hyperbolic equations,}
Mem. Amer. Math. Soc. 275 (1983).
\medskip\noindent
[M.2], A. Majda,
{\it The existence of multi-dimensional shock fronts,}
Mem. Amer. Math. Soc. 281 (1983).
\medskip\noindent
[M.3] A. Majda,
{\it Compressible fluid flow and systems of conservation laws in several
space variables,} Springer-Verlag, New York (1984), viii+ 159 pp.
\medskip\noindent
[M.4] A. Majda,
{\it A qualitative model for dynamic combustion,} 
SIAM J. Appl. Math. 41 (1981) 70--93.
\medskip\noindent
[MR] Majda and Rosales, 
{\it Weakly nonlinear detonation waves,}
SIAM J. Appl. Math. 43 (1983) 1086--1118. 
\medskip\noindent
[MZ.1] C. Mascia and K. Zumbrun,
{\it Pointwise Green's function bounds and stability of relaxation
shocks,} to appear, Indiana Math. J. (2002).
\medskip\noindent
[MZ.2] C. Mascia and K. Zumbrun,
{\it Stability of viscous shock profiles for
symmetrizable hyperbolic--parabolic systems,}
preprint (2001), available: math.indiana.edu/home/kzumbrun.
\medskip\noindent
[MZ.3] C. Mascia and K. Zumbrun,
{\it Pointwise Green's function bounds for shock profiles
with degenerate viscosity,} to appear, Arch. for Rat. Mech.
and Anal. (2003).
\medskip\noindent
[MT] U.B. McVey, U.B. and T.Y. Toong, 
{\it Mechanism of instabilities in
exothermic blunt-body flows,} Combus. Sci. Tech. 3 (1971) 63--76.
\medskip\noindent
[Me] R. Menikoff, 
{\it Determining curvature effect on detonation
velocity from rate stick experiment,}
Impact of Comp. in Sci. and Eng. 1 (1989) 168--179.
\medskip\noindent
[MeLB] 
R. Menikoff, K.S Lackner, and B.G. Bukiet, 
{\it Modeling flows with curved detonation waves,} 
Comb. and Flame 104 (1996) 219--240.
\medskip\noindent
[MeP] R. Menikoff and B. Plohr, 
{\it The Riemann problem for fluid flow of real materials,}
Rev. Modern Phys. 61 (1989), no. 1, 75--130.
\medskip\noindent
[Pa] A. Pazy, {\it Semigroups of linear operators and applications 
to partial differential equations,} Applied Mathematical Sciences, 44, 
Springer-Verlag, New York-Berlin, (1983) viii+279 pp. ISBN: 0-387-90845-5.
\medskip\noindent
[PZ] R. Plaza and K. Zumbrun,
{\it An Evans function approach to spectral stability
of small-amplitude viscous shock profiles,}
preprint (2002).
\medskip\noindent
[RV] J.-M. Roquejoffre and J.-P. Vila,
{\it Stability of ZND detonation waves in the 
Majda combustion model,} Asymptot. Anal. 18 (1998), no. 3-4, 329--348.
\medskip\noindent
[SS] Sandstede, Scheel, 
{\it Evans function and blow-up methods in critical eigenvalue problems,}
preprint (2002).
\medskip\noindent
[Sat] D. Sattinger,
{\it On the stability of waves of nonlinear parabolic systems},
Adv. Math. 22 (1976) 312--355.
\medskip\noindent
[ScSh] S. Schecter and M. Shearer, 
{\it Undercompressive shocks for non-strictly 
hyperbolic conservation laws,}
J. Dynamics Differential Equations 3 (1991), no. 2, 199--271.
\medskip\noindent
[S] M. Short,
{\it Multidimensional linear stability of a detonation wave at 
high activation energy,} SIAM J. Appl. Math. 57 (1997) 307--326.
\medskip\noindent
[Sz] A. Szepessy,
{\it Dynamics and stability of a weak detonation wave,}
Comm. Math. Phys. 202 (1999), no. 3, 547--569.
\medskip\noindent
[ZH] K. Zumbrun and P. Howard,
{\it Pointwise semigroup methods and stability of viscous shock waves,}
Indiana Mathematics Journal V47 (1998) no. 4, 741--871;
{\it Errata,} to appear, Indiana Mathematics Journal 
(2002) available math.indiana.edu/home/kzumbrun. 
\medskip\noindent
[Z.1] K. Zumbrun,  {\it Stability of viscous shock waves},
Lecture Notes, Indiana University (1998).
\medskip\noindent
[Z.2] K. Zumbrun, {\it Refined wave--tracking and nonlinear 
stability of viscous lax shocks,} Methods Appl. Anal. 7 (2000) 747--768. 
\medskip\noindent
[Z.3] K. Zumbrun, 
{\it Multidimensional stability of planar viscous shock waves,}
Advances in the theory of shock waves, 307--516, 
Progr. Nonlinear Differential Equations Appl., 47, Birkhäuser Boston,
Boston, MA, 2001.
\medskip\noindent
[Z.4] K. Zumbrun,  
{\it Multidimensional stability of Navier--Stokes shock profiles,}
Handbook of Fluid dynamics, in preparation.
\medskip\noindent
[Z.5] K. Zumbrun, 
{\it Stability index for relaxation and real viscosity
systems,} unpublished note (revised appendix for ref. [Z.3]),
(2002) available math.indiana.edu/home/kzumbrun. 
\medskip\noindent
[Z.6] K. Zumbrun, {\it Dynamical Stability of Phase Transitions 
in the P-System with Viscosity-Capillarity}, 
SIAM J. Appl. Math. 60 (2000) 1913--1924. 
\medskip\noindent
[ZPM] K. Zumbrun, B. Plohr, and D. Marchesin,
{\it Scattering behavior of transitional shock waves,}
Second Workshop on Partial Differential Equations 
(Rio de Janeiro, 1991). Mat. Contemp. 3 (1992), 191--209.
\medskip\noindent
[ZS] K. Zumbrun and D. Serre,
{\it Viscous and inviscid stability of multidimensional 
planar shock fronts,} Indiana Univ. Math. J. 48 (1999), no. 3,
937--992.
\endRefs
\enddocument